\documentclass[12pt,oneside,reqno]{amsart}
\usepackage{txfonts}
\usepackage{bbm}
\usepackage{amsmath}
\usepackage{graphicx}
\usepackage{mathrsfs}
\usepackage{stmaryrd}
\usepackage{amsfonts}
\usepackage{enumerate,amsmath,amssymb,amsthm}
\numberwithin{equation}{section}

\newcommand{\be}{\begin{eqnarray}}
\newcommand{\ee}{\end{eqnarray}}
\newcommand{\ce}{\begin{eqnarray*}}
\newcommand{\de}{\end{eqnarray*}}
\newtheorem{theorem}{Theorem}[section]
\newtheorem{lemma}[theorem]{Lemma}
\newtheorem{remark}[theorem]{Remark}
\newtheorem{definition}[theorem]{Definition}
\newtheorem{proposition}[theorem]{Proposition}
\newtheorem{Examples}[theorem]{Example}
\newtheorem{corollary}[theorem]{Corollary}

\def\eps{\varepsilon}

\def\a{\alpha}

\def\p{\partial}
\def\d{\delta}

\def\[{{\Big[}}
\def\]{{\Big]}}
\def\<{{\langle}}
\def\>{{\rangle}}
\def\({{\Big(}}
\def\){{\Big)}}

\def\bx{{\mathbf{x}}}

\def\dif{{\mathord{{\rm d}}}}

\def\no{\nonumber}
\def\={&\!\!=\!\!&}
\def\bt{\begin{theorem}}
\def\et{\end{theorem}}
\def\bl{\begin{lemma}}
\def\el{\end{lemma}}
\def\br{\begin{remark}}
\def\er{\end{remark}}
\def\bd{\begin{definition}}
\def\ed{\end{definition}}
\def\bp{\begin{proposition}}
\def\ep{\end{proposition}}
\def\bc{\begin{corollary}}
\def\ec{\end{corollary}}
\def\bx{\begin{Examples}}
\def\ex{\end{Examples}}

\def\cB{{\mathcal B}}

\def\cF{{\mathcal F}}

\def\cL{{\mathcal L}}
\def\cM{{\mathcal M}}

\def\cR{{\mathcal R}}
\def\cS{{\mathcal S}}

\def\mB{{\mathbb B}}

\def\mE{{\mathbb E}}

\def\mN{{\mathbb N}}

\def\mR{{\mathbb R}}
\def\mS{{\mathbb S}}

\def\mX{{\mathbb X}}

\def\mZ{{\mathbb Z}}

\def\sA{{\mathscr A}}

\def\sD{{\mathscr D}}

\def\sI{{\mathscr I}}

\def\sL{{\mathscr L}}
\def\sM{{\mathscr M}}

\def\sR{{\mathscr R}}

\def\geq{\geqslant}
\def\leq{\leqslant}

\def\div{\mathord{{\rm div}}}

\allowdisplaybreaks

\begin{document}

\title{Stochastic Partial Differential Equations with
Unbounded and Degenerate Coefficients}
\date{}
\author{Xicheng Zhang}

\thanks{{\it Keywords: }DiPerna-Lions theory, Stochastic partial differential equation,
Maximal principle, Nonlinear filtering}

\dedicatory{
Department of Mathematics,
Huazhong University of Science and Technology\\
Wuhan, Hubei 430074, P.R.China,\\
School of Mathematics and Statistics\\
The University of New South Wales, Sydney, 2052, Australia\\
Email: XichengZhang@gmail.com
 }

\begin{abstract}

In this article, using DiPerna-Lions theory \cite{Di-Li},
we investigate linear second order stochastic partial differential equations
with unbounded and degenerate non-smooth coefficients, and obtain several conditions for existence and
uniqueness. Moreover, we also prove the $L^1$-integrability and
a general maximal principle for generalized solutions of SPDEs.
As applications, we study nonlinear filtering problem and
also obtain the existence and uniqueness of generalized solutions for
a degenerate  nonlinear SPDE.

\end{abstract}

\maketitle

\section{Introduction}

Consider the following second order linear stochastic partial differential equation (SPDE) in $\mR^d$:
\be
\dif u_t=(\sL_tu_t+f_t)\dif t+(\sM^l_tu_t+g^l_t)\dif B^l_t, \ \
u_0(\omega,x)=\varphi(\omega,x),\label{SPDE}
\ee
where $\{B^l_t, t\geq 0\}_{l\in\mN}$ is a sequence of independent standard
Brownian motions defined on a filtered
probability space $(\Omega,\cF, P; (\cF_t)_{t\geq 0})$, and
the random partial differential operators $\sL_t(\omega)$ and $\sM^l_t(\omega)$
are given by
\be
\sL_t(\omega)u:=\p_i(a^{ij}_t(\omega,x)\p_ju)
+\p_i(b^i_t(\omega,x)u)+c_t(\omega,x)u,\label{Ad3}
\ee
where $a^{ij}=a^{ji}$ is symmetric, and
\be
\sM^l_t(\omega)u:=\sigma^{il}_t(\omega,x)\p_iu+h^l_t(x,\omega)u.\label{Ad4}
\ee
Throughout this paper, we use the following convention: when the indices $i,j,k,l$
appear twice in a product, it will be summed. Moreover, $i,j,k$ runs from $1$ to $d$ and
$l$ runs from $1$ to $\infty$. For instance,
$$
\p_i(a^{ij}\p_j u):=\sum_{i,j=1}^d\p_i(a^{ij}\p_j u),
\ \p_ib^i:=\sum_{i=1}^d\p_ib^i,\
|\sigma^{il}\xi_i|^2:=\sum_{l=1}^\infty
\Big|\sum_{i=1}^d\sigma^{il}\xi_i\Big|^2.
$$
Important notice: if we write $|\xi_i|^2$, without confusions, it always means that $\sum_{i}|\xi_i|^2$ as above.
Below, we assume that the following parabolic condition holds:
for all  $(t,\omega,x)\in [0,T]\times\Omega\times\mR^d$ and $\xi\in\mR^d$,
\be
\sA_{a,\sigma}(\xi):=2a^{ij}_t(\omega,x)\xi_i\xi_j
-|\sigma^{il}_t(\omega,x)\xi_i|^2\geq \kappa(x) |\xi_i|^2,\label{DD}
\ee
where $\kappa(x)\geq 0$ is a non-negative measurable function. If $\kappa(x)\geq \kappa_0>0$, we say
that the super-parabolic condition holds.

Let $\cM$  be the progressive $\sigma$-field on $[0,T]\times\Omega$. Let $l^2$ be
the usual Hilbert space of all sequences of square summable real numbers.
All the coefficients are always assumed to be $\cM\times\cB(\mR^d)$-measurable.
It is well known (cf. \cite[p.131, Theorem 1]{Ro}) that under super-parabolic condition, if
$a,b,\div b, c,\sigma, h$ are {\it bounded} and $\cM\times\cB(\mR^d)$-measurable functions, and
$$
f\in L^2([0,T]\times\Omega,\cM;W^{-1,2}(\mR^d)), \ \
g\in L^2([0,T]\times\Omega,\cM;L^2(\mR^d;l^2)),
$$
where $W^{m,p}(\mR^d), m\in\mZ, p>1$ denotes the usual Sobolev space,
then for any $\varphi\in L^2(\Omega,\cF_0;L^2(\mR^d))$,
there exists a unique generalized solution to SPDE (\ref{SPDE}) in the class
$$
\mX:=L^2(\Omega;C([0,T];L^2(\mR^d)))\cap L^2([0,T]\times\Omega,\cM;W^{1,2}(\mR^d)).
$$
On the other hand, in the case of $\kappa(x)\equiv 0$ (i.e, degenerate case),
if $a, b, c,\sigma, h$ are bounded and have bounded continuous derivatives up to second order, and
$$
f\in L^2([0,T]\times\Omega,\cM;W^{1,2}(\mR^d)),\ \
g\in L^2([0,T]\times\Omega,\cM;W^{2,2}(\mR^d;l^2)),
$$
then for any $\varphi\in L^2(\Omega,\cF_0;W^{1,2}(\mR^d))$, there exists a unique
generalized solution to SPDE (\ref{SPDE}) in the same class
$\mX$ (cf. \cite[p.155, Theorem 1]{Ro}).
Moreover, in the case of super-parabolic, an analytic $L^p$-theory has been established
by Krylov \cite{Kr}. But, still the {\it boundedness} assumptions on the coefficients are required.

However, the assumptions of boundedness and non-degeneracy
would become quite restrictive in applications.
For example, in nonlinear filtering, one often meets some unbounded and degenerate coefficients.
On the other hand, in the degenerate case, for solving SPDE (\ref{SPDE}), one usually
needs to assume that the coefficients are at least twice continuously differentiable as said above.
It is natural to ask whether we can remove or weaken these restrictive assumptions.
An obvious difficulty is that when $a$ is unbounded,  it is not any more true that:
$$
W^{1,2}(\mR^d)\ni u\mapsto\p_i(a^{ij}\p_j u)\in W^{-1,2}(\mR^d).
$$
Moreover, in the degenerate case, it is not expected to have any
a priori estimate for the first order
derivative of $u$ with respect to the spatial variable if the coefficients are not smooth.

Recently, Le Bris and Lions \cite{Le-Li} studied the existence and uniqueness of
deterministic Fokker-Planck equations with degenerate and irregular coefficients.
Therein, the consideration of degeneracy is motivated by the pathwise uniqueness
of SDEs with irregular coefficients and some modelling equations
in polymeric fluids. The main tool of their proofs is the DiPerna-Lions theory (cf. \cite{Di-Li})
of renormalized solutions to linear transport equations.
The aim of the present paper is to relax the assumptions on $a, b, c$
by using the DiPerna-Lions theory (cf. \cite{Di-Li}).

We mention that a general maximal principle for SPDEs has been obtained
by Krylov \cite{Kr1} under boundedness
assumptions on coefficients. A historical remark about the maximal principle of SPDEs
is also referred to \cite{Kr1}. Moreover, in \cite{Kr2}, Krylov  studied the unique
solvability of SPDE (\ref{SPDE}) with unbounded $b,c$ and bounded $a,\sigma, h$
under super-parabolic assumption. Some other well known results about SPDEs with
unbounded coefficients in weight spaces can be found in the references of \cite{Kr2}.

This paper is organized as follows:
In Section 2, we state our main results about the well-posedness of SPDE (\ref{SPDE})
under different assumptions. In Section 3, under less conditions on the coefficients,
we first prove the existence of generalized solutions. In Section 4,
we prove a general maximum principle for the generalized solutions of SPDE (\ref{SPDE})
with $g^l\equiv 0$,
which in particular implies the uniqueness of generalized solutions. Here, a commutation lemma
of DiPerna-Lions about the mollifiers plays a crucial role. In Section 5, we study the $L^1$-integrability
and weak continuity of generalized solutions constructed in Section 3.
In Section 6, we apply our results to the linear filtering equations.
In Section 7, we  prove the existence and uniqueness of
generalized solutions for a degenerate nonlinear SPDE. In the appendix, the commutation lemma of
DiPerna and Lions is proved for the reader's convenience.

\section{Statements of Main Results}

Let $W^{m,p}(\mR^d)$ be the usual real valued Sobolev space,
$W^{m,p}(\mR^d;l^2)$ the $l^2$-valued Sobolev spaces.
Let $W^{m,p}_{loc}(\mR^d)$ and $W^{m,p}_{loc}(\mR^d;l^2)$ be
the corresponding local Sobolev space. We denote by $C^\infty_0(\mR^d)$ the set of
all smooth functions over $\mR^d$ with compact supports. For a Banach space
$(\mB,\|\cdot\|_\mB)$, by $C_w([0,T];\mB)$ we denote the space of all $\mB$-valued
bounded measurable functions on $[0,T]$ that are weakly continuous
with respect to the weak topology of $\mB$. We remark that $C_w([0,T];\mB)$ is still a Banach space
under the uniform norm.

Below, we first give the notion of generalized solutions for SPDE (\ref{SPDE}). For this, we need
to assume that
$$
\left\{
\begin{aligned}
&a^{ij},\p_ia^{ij}, b^i, c\in L^1(0,T; L^2(\Omega;L^2_{loc}(\mR^d))),\\
&\sigma^{i\cdot}, \p_i\sigma^{i\cdot}, h^\cdot\in
L^2([0,T]\times\Omega;L^2_{loc}(\mR^d; l^2)),\\
& f\in L^1([0,T]\times\Omega;L^1_{loc}(\mR^d)),\\
&g\in L^2([0,T]\times\Omega;L^1_{loc}(\mR^d; l^2)).
\end{aligned}
\right.\leqno{\bf (BasicA)}.
$$
In what follows, these assumptions will be always made if there is no special declaration,
and without confusions, we shall drop the arguments $(t,\omega,x)$. For example,
for a function $u$, we may write
$$
\int^t_0\!\!\!\int u:=\int^t_0\!\!\!\int u_s\dif x\dif s:=\int^t_0\!\!\!\int_{\mR^d} u_s(\omega,x)\dif x\dif s.
$$

\bd\label{Def1}
Let $u_0\in L^2(\Omega,\cF_0;L^2_{loc}(\mR^d))$. An $\cM\times\cB(\mR^d)$-measurable
process
$$
u\in L^\infty(0,T; L^2(\Omega;L^2_{loc}(\mR^d)))
$$
is called a generalized (or distribution) solution of SPDE (\ref{SPDE})
if for all $\phi\in C^\infty_0(\mR^d)$,
it holds that for $(\dif t\times P)$-almost all $(t,\omega)\in[0,T]\times\Omega$,
\be
\int u_t\phi\dif x\=\int u_0\phi\dif x+\int^t_0\!\!\!\int u_s\sL^*_s\phi\dif x\dif s
+\int^t_0\!\!\!\int f_s\phi\dif x\dif s\no \\
&&+\int^t_0\!\!\!\int u_s\sM^{l*}_s\phi\dif x\dif W^l_s
+\int^t_0\!\!\!\int g^l_s\phi\dif x\dif W^l_s,\label{Gen}
\ee
where $\sL^*_t$ and $\sM^{l*}_t$ are their respective adjoint operators and given by
\be
\sL^*_t(\omega)\phi:=\p_i(a^{ij}_t(\omega,x)\p_j\phi)
+b^i_t(\omega,x)\p_i\phi+c_t(\omega,x)\phi\label{Ad1}
\ee
and
\be
\sM^{l*}_t(\omega)\phi:=\p_i(\sigma^{il}_t(\omega,x)\phi)+h^l_t(\omega,x)\phi.\label{Ad2}
\ee
\ed
\br
It is easy to see that both sides of (\ref{Gen}) are well defined under the above
described basic assumptions.
\er

We now state our first result under non-degenerate assumption,
which is a direct conclusion of Propositions \ref{Th1}, \ref{Pro4},
 \ref{Pro1}, \ref{Pro7} and \ref{Pro8} below.
\bt\label{Main1}
Let parabolic condition $\sA_{a,\sigma}(\xi)\geq\kappa |\xi|^2$ be fulfilled with
$\kappa\in C^1(\mR^d;(0,\infty))$, having continuous first order derivatives.
Assume that the following conditions hold:
$$
\left\{
\begin{aligned}
&\frac{|a^{ij}|}{1+|x|^2}, \frac{|\p_j a^{ij}|}{1+|x|}
\in L^\infty([0,T]\times\Omega\times\mR^d);\\
&\frac{|b^i|}{1+|x|}, \div b\in L^2(0,T; L^\infty(\Omega\times\mR^d));\\
&c\in L^1(0,T; L^2(\Omega;L^2_{loc}(\mR^d))),\ c^+\in L^1(0,T; L^\infty(\Omega\times\mR^d));\\
&\sigma^{i\cdot}, \p_i\sigma^{i\cdot}, h, \p_k h\in
L^\infty([0,T]\times\Omega\times\mR^d;l^2),
\end{aligned}
\right.
$$
where $c^+=\max(0,c)$.
Then for any $u_0\in L^2(\Omega,\cF_0; L^2(\mR^d))$ and
$$
f\in L^2([0,T]\times\Omega\times\mR^d),
\ \ g\in L^2([0,T]\times\Omega; W^{1,2}(\mR^d;l^2)),
$$
there exists a unique generalized solution $u\in L^2(\Omega; C_w([0,T]; L^2(\mR^d)))$
to SPDE (\ref{SPDE}) in the class that
$$
\mE\left(\int^T_0\!\!\!\int\kappa |\p_iu |^2 \right)<+\infty.
$$
Moreover, in addition to the above assumptions on the coefficients and $u_0, f$,
\begin{enumerate}[{\rm (I)}]
\item if $g^l\equiv 0$, $u_0\in L^1(\Omega,\cF_0;L^1(\mR^d))$ and
$f\in L^1([0,T]\times\Omega\times\mR^d)$, then
$$
u\in L^1(\Omega; C_w([0,T]; L^1(\mR^d)))
$$
and for some $C>0$ independent of $u_0$ and $f$,
$$
\mE\left(\sup_{t\in[0,T]}\int|u_t|\right)
\leq C\mE\int|u_0|+C\mE\int^T_0\!\!\!\int |f|;
$$
\item if $f\geq 0$, $g^l\equiv 0$ and $u_0\geq 0$, then
$$
u_t(\omega,x)\geq 0,\ \ (\dif t\times P\times\dif x)-a.s.
$$
\end{enumerate}
\et

In the degenerate case, we present three different results.
The first one is a conseqeunce of Propositions \ref{Th1}, \ref{Pro6},
\ref{Pro1}, \ref{Pro7} and \ref{Pro8}.
\bt\label{Main3}
Let $a$ and $\sigma$ be independent of $x$. Assume that  the following conditions hold: for some $q>1$
$$
\left\{
\begin{aligned}
&a^{ij}\in L^1(0,T; L^\infty(\Omega));\ \
\sigma^{i\cdot}\in L^{2q}(0,T;L^\infty(\Omega; l^2));\\
&\frac{|b^i|}{1+|x|} \in L^1(0,T; L^2(\Omega\times\mR^d))
\cup L^1(0,T; L^\infty(\Omega\times\mR^d));\\
&\p_k b^i, c\in L^1(0,T; L^2(\Omega; L^{2}_{loc}(\mR^d)));
\div b, c^+\in L^1(0,T; L^\infty(\Omega\times\mR^d));\\
&h,\p_k h\in L^{2q}(0,T; L^\infty(\Omega\times\mR^d;l^2)).
\end{aligned}
\right.
$$
Then for any $u_0\in L^2(\Omega,\cF_0; L^2(\mR^d))$ and
$$
f\in L^2([0,T]\times\Omega\times\mR^d),
\ \ g\in L^2([0,T]\times\Omega; W^{1,2}(\mR^d;l^2)),
$$
there exists a unique generalized solution $u\in L^2(\Omega; C_w([0,T]; L^2(\mR^d)))$
to SPDE (\ref{SPDE}). Moreover, the same conclusions (I) and (II)
in Theorem \ref{Main1} still hold.
\et
The following result is an extension of Krylov and Rozovskii's result \cite{Kr-Ro1} (see
\cite[p.155, Theorem 1]{Ro}),
which is a consequence of Propositions \ref{Pro2}, \ref{Pro4}, \ref{Pro1}, \ref{Pro7} and \ref{Pro8}.
\bt\label{Main4}
Assume that for some $C_0>0$ and $q>1$
$$
\left\{
\begin{aligned}
&\p_k\p_ka^{ij}\xi_i\xi_j\leq C_0|\xi|^2;\ \frac{|a^{ij}|}{1+|x|^2},
\frac{|\p_k a^{ij}|}{1+|x|},\frac{|b^i|}{1+|x|},\\
&\p_k b,\p_k\div b, c,\p_k c\in L^1(0,T; L^\infty(\Omega\times\mR^d));\\
&\sigma^{i\cdot}, \p_k\sigma^{i\cdot}, \p_k\p_j\sigma^{i\cdot},h,
\p_k h, \p_k\p_jh\in L^{2q}(0,T; L^\infty(\Omega\times\mR^d;l^2)).
\end{aligned}
\right.
$$
Then for any $u_0\in L^2(\Omega,\cF_0; W^{1,2}(\mR^d))$ and
$$
f\in L^2([0,T]\times\Omega; W^{1,2}(\mR^d)), \ \ g\in L^2([0,T]\times\Omega; W^{2,2}(\mR^d;l^2)),
$$
there exists a unique generalized solution $u\in L^2(\Omega; C_w([0,T];W^{1,2}(\mR^d)))$
to SPDE (\ref{SPDE}). Moreover, the conclusions (I) and (II) in Theorem \ref{Main1} still hold.
\et

The following result is an easy consequence of
Propositions \ref{Pro3}, \ref{Pro5}, \ref{Pro1}, \ref{Pro7} and \ref{Pro8}.
\bt\label{Main2}
 Let $a^{ij}$ be given as follows
$$
a^{ij}_t(\omega,x)=\hat\sigma^{il}_t(\omega,x)\hat\sigma^{jl}_t(\omega,x)
$$
such that for some $\a>1/2$
$$
|\hat\sigma^{il}_t\xi_i|^2\geq \a\cdot|\sigma^{il}_t\xi_i|^2.
$$
Assume that the following conditions hold: for some $q>1$
$$
\left\{
\begin{aligned}
&\frac{\hat\sigma^{i\cdot}}{1+|x|}, \p_i\hat\sigma^{i\cdot}
\in L^2(0,T; L^\infty(\Omega\times\mR^d; l^2));\\
&\frac{|b^i|}{1+|x|}, \div b, c^+\in L^1(0,T; L^\infty(\Omega\times\mR^d));\\
&\p_k b^i, c\in L^1(0,T; L^2(\Omega;L^{2}_{loc}(\mR^d)));\\
&\sigma^{i\cdot}, \p_i\sigma^{i\cdot}, h
\in L^{2q}(0,T; L^\infty(\Omega\times\mR^d; l^2)).
\end{aligned}
\right.
$$
Then for any $u_0\in L^2(\Omega,\cF_0; L^2(\mR^d))$ and
$$
f\in L^2([0,T]\times\Omega\times\mR^d),
\ \ g\in L^2([0,T]\times\Omega\times\mR^d; l^2),
$$
there exists a unique generalized solution $u\in L^2(\Omega; C_w([0,T]; L^2(\mR^d)))$
to SPDE (\ref{SPDE}) satisfying
\be
&&\mE\left(\sup_{s\in[0,T]}\int|u_{s}|^{2}\right)
+\mE\left(\int^T_0\!\!\!\int\big(|\hat\sigma^{il}\p_i u|^2+|\sigma^{il}\p_i u|^2\big)\right)\no\\
&&\qquad\qquad\leq C\left(\mE\int|u_0|^2+\mE\int^T_0\!\!\!\int|f|^2
+\mE\int^T_0\!\!\!\int|g|^2\right),\label{PL1}
\ee
where $C$ is independent of $u_0, f$ and $g$.
Moreover, the conclusions (I) and (II) in Theorem \ref{Main1} still hold.
\et

\section{Existence of Generalized Solutions}

In the sequel, we shall use the following conventions:
The letter $C$ denotes a constant whose value may change in different occasions,
and $\ell_t$ denotes an $L^1$-integrable real function on $[0,T]$
which may be different in different lines.

We now state our first existence result.
\bp\label{Th1}
Assume that {\bf (BasicA)} and the following conditions hold: for some $q>1$
\be
&a^{ij}\in L^1(0,T; L^\infty(\Omega;L^2_{loc}(\mR^d))),\\
&\div b,\ c^+\in L^1(0,T; L^\infty(\Omega\times\mR^d)),\label{Eq2}\\
&\sigma^{i\cdot}, \p_i\sigma^{i\cdot}, h, \p_k h\in L^{2q}(0,T; L^\infty(\Omega\times\mR^d;l^2)).\label{Eq3}
\ee
Then for any $u_0\in L^{2}(\Omega,\cF_0;L^{2}(\mR^d))$ and
$$
f\in L^2([0,T]\times\Omega\times\mR^d),
\ \ g\in L^2([0,T]\times\Omega; W^{1,2}(\mR^d;l^2)),
$$
there exists a generalized  solution
$$u\in L^2(\Omega; L^\infty(0,T;L^2(\mR^d)))$$ to SPDE (\ref{SPDE}).
Moreover, if $\kappa\in C^1(\mR^d;(0,\infty))$, then the above
generalized solution  also satisfies
\be
\mE\left(\int^T_0\!\!\!\int\kappa |\p_i u|^2\right)<+\infty.\label{Es6}
\ee
\ep
\br
If $\kappa(x)\geq\kappa_0>0$, then the above assumptions can be weakened.
Since there is a complete theory in the super-parabolic case (cf. \cite{Ro,Kr}),
it is not pursued here.
\er

For proving this proposition, we adopt the argument of mollifying the coefficients.
Let $\rho\in C^\infty_0(\mR^d)$ be a regularizing kernel function with
$$
\mathrm{supp}(\rho)\subset \bar B,\ \ \rho>0\mbox{\ on\ $B_1$},\ \ \int\rho=1,
$$
where $B_1:=\{x\in\mR^d: |x|<1\}$.
Let $\chi\in C^\infty_0(\mR^d)$ be a non-negative cutoff function
with $\chi=1$ on the unit ball and $\chi=0$ outside the ball of radius $2$.
Set for $\eps\in(0,1)$
$$
\rho_{\eps}(x):=\eps^{-d}\rho(\eps^{-1}x),\ \ \chi_\eps:=\chi(\eps x).
$$
Define
\be\label{de2}
\left\{
\begin{aligned}
&a^{ij}_{t,\eps}:=(a^{ij}_t*\rho_{\eps}) \chi^2_\eps,&
\sigma^{il}_{t,\eps}:=(\sigma^{il}_t*\rho_{\eps})  \chi_\eps,\\
&c_{t,\eps}:=(c_t*\rho_{\eps})) \chi_\eps,&
h^l_{t,\eps}:=(h^l_t*\rho_{\eps}) \chi_\eps
\end{aligned}
\right.
\ee
and
$$
f_{t,\eps}:=(f_t*\rho_{\eps}) \chi_\eps,\ \
g^l_{t,\eps}:=(g^l_t*\rho_{\eps}) \chi_\eps,
$$
where the asterisk stands for the convolution in $x$.
Moreover, we  define
\be
b^i_{t,\eps}:=[(b^i_t\wedge(1/\eps))\vee(-1/\eps)]*\rho_{\eps}.\label{de1}
\ee
\br
Here, for a vector field $b$, we directly truncate $b$ rather than
multiplying a cutoff function on $\mR^d$ so that
$\|\mathrm{div}b_{t,\eps}\|_\infty\leq\|\mathrm{div}b_t\|_\infty$. Otherwise, we need
an extra assumption on $b$ (see (\ref{Es2}) and (\ref{PP1}) below).
\er

We need the following simple lemma.
\bl\label{Le1}
(i) Let parabolic condition $\sA_{a,\sigma}(\xi)\geq\kappa |\xi|^2$ hold.
Set $\kappa_\eps:=\kappa*\rho_\eps$. Then,
\be
\sA_{a_\eps,\sigma_\eps}(\xi)\geq \kappa_\eps \chi_\eps^2
|\xi_i|^2,\ \forall \xi\in\mR^d,\label{PL2}
\ee
where $\sA_{a_\eps, b_\eps}$ is defined by (\ref{DD}).

(ii) Assume that (\ref{Eq2}) and (\ref{Eq3}) hold. Then, for some
$\ell_t\in L^1(0,T)$, \be
\sup_{\eps\in(0,1)}\||\p_ib^i_{t,\eps}|+|c^+_{t,\eps}|\|_{L^\infty(\Omega\times\mR^d)}\leq
\ell_t\label{Es2} \ee and \be
\sup_{\eps\in(0,1)}\||\sigma^{i\cdot}_{t,\eps}|+|\p_i\sigma^{i\cdot}_{t,\eps}|
+|h_{t,\eps}|+|\nabla
h_{t,\eps}|\|_{L^\infty(\Omega\times\mR^d;l^2)}^{2}\leq\ell_t.\label{Es3}
\ee

(iii) Let $\sL^*_{t,\eps}$ and $\sM^{l*}_{t,\eps}$ be defined
in terms of $a_{t,\eps},b_{t,\eps},c_{t,\eps}$
and $\sigma_{t,\eps},h_{t,\eps}$
as in (\ref{Ad1}) and (\ref{Ad2}). Then for any $\phi\in C^\infty_0(\mR^d)$,
\be
\lim_{\eps\to 0}\int^T_0\|\sL^*_{t,\eps}\phi-\sL^*_t\phi\|_{L^2(\Omega\times\mR^d)}\dif t=0\label{Lim1}
\ee
and
\be
\lim_{\eps\to 0}\int^T_0\|\sM^{l*}_{t,\eps}\phi-\sM^{l*}_t\phi\|^2_{L^2(\Omega\times\mR^d)}\dif t=0.\label{Lim2}
\ee
\el
\begin{proof}
(i). By virtue of $\int\rho_{\eps}=1$, we have
\be
|\sigma^{il}_{t,\eps}\xi_i|^2=|\sigma^{il}_t\xi_i* \rho_{\eps}|^2 \chi_\eps^2
\leq(|\sigma^{il}_t\xi_i|^2*\rho_{\eps}) \chi_\eps^2.\label{Es9}
\ee
Hence,
\ce
\sA_{a_\eps,\sigma_\eps}(\xi)&=&2a^{ij}_{t,\eps}\xi_i\xi_j-|\sigma^{il}_{t,\eps}\xi_i|^2\\
&\geq&\left(2a^{ij}_t\xi_i\xi_j-|\sigma^{il}_t\xi_i|^2\right)*\rho_{\eps}\cdot\chi_\eps^2\\
&\geq& ((\kappa |\xi_i|^2)*\rho_{\eps})\chi_\eps^2=\kappa_\eps\chi_\eps^2|\xi_i|^2.
\de

(ii). Estimate (\ref{Es2}) is direct from  definition (\ref{de1}) and (\ref{Eq2}).
Estimate (\ref{Es3}) follows from
\be
|\p_i\chi_\eps(x)|\leq C\eps  1_{[1/\eps,2/\eps]}(|x|)\label{PL5}
\ee
and (\ref{Eq3}).

(iii). Limits (\ref{Lim1}) and (\ref{Lim2}) follow from the property of convolution mollifying.
\end{proof}

Consider now the following approximation equation:
\be
\dif u_{\eps,t}=(\sL_{t,\eps} u_{\eps,t}+f_{t,\eps})\dif t
+(\sM^l_{t,\eps} u_{\eps,t}+g^l_{t,\eps})\dif B^l_t, \label{App}
\ee
subject to $u_{\eps,0}:=(u_0*\rho_\eps)\chi_\eps$, where
$\sL_{t,\eps}$ and $\sM^{l}_{t,\eps}$ are defined respectively in terms of $a_{t,\eps},b_{t,\eps},c_{t,\eps}$
and $\sigma_{t,\eps},h_{t,\eps}$ as in (\ref{Ad3}) and (\ref{Ad4}).
Notice that all the coefficients of (\ref{App}) are  smooth in $x$,
and their derivatives of all orders in $x$ are uniformly bounded in $(\omega,x)$ for fixed $t$.
In fact, we may further assume that all the coefficients together with all of their
derivatives in $x$ are uniformly bounded in $(t,\omega,x)$
if we also mollify the time variable and cut off it as done for $x$. We omit this tedious step
for the sake of simplicity.
Moreover, if we let $W^\infty(\mR^d)=\cap_{k\in\mN}W^{k,2}(\mR^d)$, then $f_\eps,g^l_\eps\in
L^2([0,T]\times\Omega;W^\infty(\mR^d))$ and $u_{\eps,0}\in L^2(\Omega,\cF_0; W^\infty)$.
Thus, by \cite[p.155, Theorem 1]{Ro}, there exists a unique smooth solution
$u_\eps\in L^2(\Omega;C([0,T];W^\infty(\mR^d)))$ to equation (\ref{App}).

Below, for the simplicity, we sometimes drop the time variable $t$ in
$a_{t,\eps}, b_{t,\eps}$, etc.
By It\^o's formula and the integration by parts formula, we have
\be
\dif\int u_{\eps}^2&=&
\left(2\int u_\eps (\sL_\eps u_\eps+f_\eps)+\int|\sM^l_\eps u_\eps+g^l_\eps|^2\right)\dif t\no\\
&&+2\left(\int u_\eps (\sM^l_\eps u_\eps+g^l_\eps)\right)\dif B^l_t\no\\
&=&\left(-\int\sA_{a_\eps,\sigma_\eps}(\nabla u_\eps)
+2\int u_\eps (\p_i(b^i_\eps u_\eps)+c_\eps u_\eps+f_\eps)\right)\dif t\no\\
&&+\left(\int\Big(2\sigma^{il}_\eps\p_iu_\eps
 (h^l_\eps u_\eps+g^l_\eps)+(h^l_\eps u_\eps+g^l_\eps)^2\Big)\right)\dif t\no\\
&&+\left(\int\Big(u_\eps^2\p_i\sigma^{il}_\eps+2u_\eps
(h^l_\eps u_\eps+g^l_\eps)\Big)\right)\dif B^l_t.\label{Es01}
\ee
Observing that
$$
2\int u_\eps\p_i(b^i_\eps u_\eps)=\int u_\eps^2\p_ib^i_\eps
$$
and
\be
2\int\p_iu_\eps\sigma^{il}_\eps(h^l_\eps u_\eps+g^l_\eps)
=-\int u_\eps^2\p_i(\sigma^{il}_\eps h^l_\eps)
-2\int u_\eps\p_i(\sigma^{il}_\eps g^l_\eps),\label{L6}
\ee
by integrating both sides of (\ref{Es01}) in time from $0$ to $t$, we get
\be
\int u_{\eps,t}^2&=&\int u_{\eps,0}^2
-\int^t_0\!\!\!\int \sA_{a_\eps,\sigma_\eps}(\nabla u_\eps)
+\int^t_0\!\!\!\int u_\eps^2\p_i(b^i_\eps-\sigma^{il}_\eps h^l_\eps),\no\\
&&+\int^t_0\!\!\!\int\Big(2u_\eps (c_\eps u_\eps+f_\eps-\p_i(\sigma^{il}_\eps g^l_\eps))+
(h^l_\eps u_\eps+g^l_\eps)^2\Big)\no\\
&&+\int^t_0\left(\int\big(u_\eps^2\p_i\sigma^{il}_\eps+2u_\eps
(h^l_\eps u_\eps+g^l_\eps)\big)\right)\dif B^l_s.\label{Es1}
\ee

We are now in a position to give:

\vspace{3mm}

{\it Proof of Proposition \ref{Th1}:}
By (\ref{Es1}) and Lemma \ref{Le1}, we have
\ce
\int|u_{\eps,t}|^{2}&\leq&\int|u_{\eps,0}|^{2}-\int^t_0\!\!\!\int\chi_\eps^2\kappa_\eps
|\p_i u_\eps|^2+\int^t_0\ell_s\left(1+\int|u_\eps|^{2}\right)\dif s\\
&&+\int^t_0\left(\int\big(u_\eps^2\p_i\sigma^{il}_\eps+2u_\eps
(h^l_\eps u_\eps+g^l_\eps)\big)\right)\dif B^l_s.
\de
First taking supremum in time and then expectations, by Burkholder's inequality, we get
\ce
&&\mE\left(\sup_{s\in[0,t]}\int|u_{\eps,s}|^{2}\right)
+\mE\int^t_0\!\!\!\int\chi_\eps^2\kappa_\eps|\p_i u_{\eps,s}|^2\leq \\
&&\qquad\qquad\leq\int|u_{\eps,0}|^2+\int^t_0\ell_s\left(1+\mE\int|u_{\eps,s}|^{2}\right)\dif s\\
&&\qquad\qquad\quad+C\mE\left(\int^t_0\left|\int\Big(u_\eps^2\p_i\sigma^{il}_\eps+2u_\eps
(h^l_\eps u_\eps+g^l_\eps)\Big)\right|^2\dif s\right)^{1/2}.
\de
Here and below, the constant $C$ is independent of $\eps$.
The last term denoted by $\sI$ can be controlled as follows: by (\ref{Es3}) and Young's inequality
\ce
\sI&\leq&C\mE\left(\int^t_0\left(\int|u_{\eps,s}|^2\right)
\left(\ell_s\int|u_{\eps,s}|^2+\int |g^l_{\eps,s}|^2\right)\dif s\right)^{1/2}\\
&\leq&C\mE\left(\sup_{s\in[0,t]}\int|u_{\eps,s}|^2
\int^t_0\left(\ell_s\int|u_{\eps,s}|^2+\int |g^l_{s}|^2\right)\dif s\right)^{1/2}\\
&\leq&\frac{1}{2}\mE\left(\sup_{s\in[0,t]}\int|u_{\eps,s}|^2\right)
+C\int^t_0\ell_s \mE\left(\int|u_{\eps,s}|^2\right)\dif s+C\mE\int^T_0\!\!\!\int |g^l|^2.
\de
 Combining the above calculations, we obtain
\ce
&&\mE\left(\sup_{s\in[0,t]}\int|u_{\eps,s}|^{2}\right)
+\mE\int^t_0\!\!\!\int\chi_\eps^2\kappa_\eps|\p_i u_{\eps,s}|^2\\
&&\qquad\leq C+C\int^t_0\ell_s\mE\left(\int|u_{\eps,s}|^2\right)\dif s\\
&&\qquad\leq C+C\int^t_0\ell_s\mE\left(\sup_{r\in[0,s]}\int|u_{\eps,r}|^2\right)\dif s.
\de
By Gronwall's inequality,
\be
\mE\left(\sup_{s\in[0,T]}\int|u_{\eps,s}|^{2}\right)
+\mE\int^T_0\!\!\!\int\chi_\eps^2\kappa_\eps
|\p_i u_{\eps,s}|^2\leq C.\label{Es5}
\ee

Consider now the Banach space $\mB_1:=L^\infty(0,T;L^2(\Omega\times\mR^d)))$
and the reflexive Banach space $\mB_2:=L^{2p}(0,T;L^2(\Omega\times\mR^d))$,
where $p=\frac{q}{q-1}\in(1,\infty)$. The sequence
$u_\eps$ is then uniformly bounded in $\mB_1\subset\mB_2$. So, there exists a $u\in\mB_1$ and a subsequence
$u_{\eps_k}$ such that $u_{\eps_k}$ weakly *  in $\mB_1$ (weakly in $\mB_2$)
converges to $u$.
Let $\phi\in C^\infty_0(\mR^d)$ and $\ell\in L^\infty([0,T]\times\Omega)$. Then by (\ref{App}), we have
\ce
\mE\int^T_0\!\!\!\int u_{\eps,t}\phi\ell_t\dif x\dif t
&=&\mE\int^T_0\!\!\!\int u_{\eps,0}\phi\ell_t\dif x\dif t+
\mE\int^T_0\ell_t\int^t_0\!\!\!\int u_{\eps,s}\sL^*_{s,\eps}\phi\dif x\dif s\dif t\no\\
&&+\mE\int^T_0\ell_t\int^t_0\!\!\!\int f_{s,\eps}\phi\dif x\dif s\dif t\no \\
&&+\mE\int^T_0\ell_t\int^t_0\!\!\!\int u_{\eps,s}\sM^{l*}_{s,\eps}\phi\dif x\dif W^l_s\dif t\no\\
&&+\mE\int^T_0\ell_t\int^t_0\!\!\!\int g^l_{s,\eps}\phi\dif x\dif W^l_s\dif t.
\de
We want to take limits $\eps\to 0$ for both sides of the above equality.
Let us first prove that
\be
\mE\int^T_0\ell_t\int^t_0\!\!\!\int u_{\eps,s}\sM^{l*}_{s,\eps}\phi\dif x\dif W^l_s\dif t
\stackrel{\eps\to 0}{\longrightarrow}\mE\int^T_0\ell_t\int^t_0\!\!\!\int u_{s}
\sM^{l*}_{s}\phi\dif x\dif W^l_s\dif t.\label{Ee1}
\ee
By (\ref{Lim2}) and (\ref{Es5}), it is easy to see that
$$
\mE\left|\int^T_0\ell_t\int^t_0\!\!\!\int u_{\eps,s}(\sM^{l*}_{s,\eps}\phi-\sM^{l*}_{s}\phi)
\dif x\dif W^l_s\dif t\right|\stackrel{\eps\to 0}{\longrightarrow}0.
$$
For $u\in\mB_2$, we define
$$
(\cR u)_t:=\int^t_0\!\!\!\int u_s\sM^{l*}_{s}\phi\dif x\dif W^l_s.
$$
By Burkholder's inequality and H\"older's inequality, we have
\ce
\mE\int^T_0|(\cR u)_t|^2\dif t&\leq&T\int^T_0\mE\left|\int u_s\sM^{l*}_{s}\phi\dif x\right|^2\dif s\\
&\leq&T\left(\int^T_0\left(\mE\int |u_s|^{2}\right)^{p}\dif s\right)^{1/p}
\left(\int^T_0\left(\mE\int|\sM^{l*}_{s}\phi|^{2}
\right)^q\dif s\right)^{1/q}\\
&\stackrel{(\ref{Eq3})}{\leq}&C_\phi T\left(\int^T_0\left(\mE\int |u_s|^{2}\right)^{p}\dif s\right)^{1/p}
=C_\phi T\|u\|_{\mB_2}^2,
\de
which means that $\cR:\mB_2\to L^2([0,T]\times\Omega)$ is a strongly continuous operator. So, $\cR$
is also weakly continuous, and
$$
\mE\int^T_0\ell_t(\cR u_\eps)_t\dif t
\stackrel{\eps\to 0}{\longrightarrow}\mE\int^T_0\ell_t (\cR u)_t\dif t.
$$
Thus, (\ref{Ee1}) is proven. By Lemma \ref{Le1} and passing to limits, as above, we finally obtain
\ce
\mE\int^T_0\!\!\!\int u_{t}\phi\ell_t\dif x\dif t
&=&\mE\int^T_0\!\!\!\int u_{0}\phi\ell_t\dif x\dif t+
\mE\int^T_0\ell_t\int^t_0\!\!\!\int u_{s}\sL^*_{s}\phi\dif x\dif s\dif t\no\\
&&+\mE\int^T_0\ell_t\int^t_0\!\!\!\int f_{s}\phi\dif x\dif s\dif t\no \\
&&+\mE\int^T_0\ell_t\int^t_0\!\!\!\int u_{s}\sM^{l*}_{s}\phi\dif x\dif W^l_s\dif t\no\\
&&+\mE\int^T_0\ell_t\int^t_0\!\!\!\int g^l_{s}\phi\dif x\dif W^l_s\dif t.
\de
Equality (\ref{Gen}) then follows by the arbitrariness of $\ell_t\in L^\infty([0,T]\times\Omega)$.

We now prove $u\in L^2(\Omega; L^\infty(0,T; L^2(\mR^d)))$.
By Banach-Saks theorem (cf. \cite{Du-Sc}),
there exists another subsequence (still denoted by $\eps_k$) such that
its Ces\`aro mean $\bar u_{\eps_n}:=\frac{\sum_{k=1}^n u_{\eps_k}}{n}$ strongly converges
to $u$ in $\mB_2$. Thus, there exist a subsequence still denoted by $\eps_n$
and a null set $A\subset[0,T]\times\Omega$ such that for all $(t,\omega)\notin A$
$$
\lim_{n\to\infty}\|\bar u_{\eps_n,t}(\omega)-u_t(\omega)\|_{L^2(\mR^d)}=0.
$$
Let $S_\omega:=\{t\in[0,T]: (t,\omega)\in A^c\}$ be the section of $A^c$. By Fubini's theorem,
for $P$-almost all $\omega$, $S_\omega$ has full Lebesgue measure.
Thus,
$$
\sup_{t\in S_\omega}\|u_{t}(\omega)\|_{L^2(\mR^d)}\leq\varliminf_{n\to\infty}\sup_{t\in S_\omega}
\|\bar u_{\eps_n,t}(\omega)\|_{L^2(\mR^d)}
\leq\varliminf_{n\to\infty}\frac{1}{n}\sum_{k=1}^n\sup_{t\in S_\omega}
\|u_{\eps_k,t}(\omega)\|_{L^2(\mR^d)},
$$
which together with (\ref{Es5}) yields
\be
\mE\left(\mathrm{ess.}\sup_{t\in[0,T]}\int |u_t|^{2}\right)<+\infty.\label{Es4}
\ee

Let $\sD\subset C^\infty_0(\mR^d;\mR^d)$ be a countable and dense subset of $L^2(\mR^d;\mR^d)$.
Noting that for fixed $\phi\in \sD$ and for $(\dif t\times P)$-almost
all $(t,\omega)\in[0,T]\times\Omega$,
\ce
\int u_t(\omega)\p_i(\sqrt{\kappa}\phi^i)&=&\lim_{n\to\infty}
\int \bar u_{\eps_n,t}(\omega)\p_i(\sqrt{\kappa}\phi^i)\\
&=&\lim_{n\to\infty}\frac{1}{n}\sum_{k=1}^n\int
u_{\eps_k,t}(\omega)\p_i(\chi_{\eps_k}\sqrt{\kappa_{\eps_k}}\phi^i),
\de
we have for $(\dif t\times P)$-almost all $(t,\omega)\in[0,T]\times\Omega$
\ce
\left(\int\kappa|\p_i u_t(\omega)|^2\right)^{1/2}&=&
\sup_{\phi\in\sD}\frac{1}{\|\phi\|_{L^2}}\int u_t(\omega)\p_i(\sqrt{\kappa}\phi^i)\\
&\leq&\varliminf_{n\to\infty}\frac{1}{n}\sum_{k=1}^n\sup_{\phi\in\sD}\frac{1}{\|\phi\|_{L^2}}
\int u_{\eps_k,t}(\omega)\p_i(\chi_{\eps_k}\sqrt{\kappa_{\eps_k}}\phi^i)\\
&=&\varliminf_{n\to\infty}\frac{1}{n}\sum_{k=1}^n\left(
\int \chi_{\eps_k}^2\kappa_{\eps_k}|\p_i u_{\eps_k,t}(\omega)|^2\right)^{1/2}.
\de
Thus, by (\ref{Es5}), we get (\ref{Es6}). The proof of Proposition \ref{Th1} is complete.
\begin{flushright} $\square$\end{flushright}


For proving the uniqueness, we need more regular solutions. Below,
we give two such results in the degenerate  case.
The first one is an extension of Krylov and Rozovskii's result \cite{Kr-Ro1}
(see also \cite[p.155, Theorem 1]{Ro}).
Therein, an Oleinik's lemma (see \cite[p.44, Lemma 2.4.3]{St}
and \cite[p.161, Proposition 3]{Ro}) plays a crucial role.
\bp\label{Pro2}
Assume that the following conditions hold: for some $q>1$
\be
\p_k\p_ka^{ij}\xi_i\xi_j\leq C_0|\xi_i|^2,\ \frac{|a^{ij}|}{1+|x|^2}, \frac{|\p_k a^{ij}|}{1+|x|}
\in L^1(0,T; L^\infty(\Omega\times\mR^d)),\label{P1}
\ee
\be
\frac{|b^i|}{1+|x|}, \p_k b,\p_k\div b, c,\p_k c
\in L^1(0,T; L^\infty(\Omega\times\mR^d)),\label{P2}
\ee
\be
\sigma^{i\cdot}, \p_k\sigma^{i\cdot}, \p_k\p_j\sigma^{i\cdot},h,
\p_k h, \p_k\p_jh\in L^{2q}(0,T; L^\infty(\Omega\times\mR^d;l^2)).\label{P3}
\ee
Then the generalized solution constructed in Proposition \ref{Th1} also satisfies
\be
\mE\left(\mathrm{ess.}\sup_{t\in[0,T]}\int|\p_i u_t|^2\right)<+\infty.\label{PL4}
\ee
\ep
\begin{proof}
By differentiating SPDE (\ref{App}) in the $k$th spatial coordinate $x_k$, we obtain
\be
\dif\p_ku_\eps&=&\Big(\sL_\eps \p_ku_\eps+[\p_k,\sL_\eps](u_\eps)+\p_kf_\eps\Big)\dif t\no\\
&&+\Big(\sM^l_\eps \p_ku_\eps+[\p_k,\sM^l_\eps](u_\eps)+\p_kg^l_\eps\Big)\dif B^l_t,  \label{App1}
\ee
where
$$
[\p_k,\sL_\eps](u)=\p_k(\sL_\eps(u))-\sL_\eps(\p_k u)=\p_i(\p_ka^{ij}_\eps\p_j u)
+\p_i(\p_k b^i_\eps u)+(\p_k c_\eps) u
$$
and
$$
[\p_k,\sM^l_\eps](u)=\p_k(\sM^l_\eps(u))-\sM^l_\eps(\p_ku)=
\p_k\sigma^{il}_\eps\p_i u+(\p_k h^l_\eps) u.
$$

Similar to (\ref{Es1}), we have
\be
\int |\p_k u_{\eps,t}|^2
&=&\int |\p_k u_{\eps,0}|^2
-\int^t_0\!\!\!\int\Big(2a^{ij}_\eps\p_i \p_ku_\eps\p_j \p_ku_\eps
-|\sigma^{il}_\eps\p_i \p_ku_\eps|^2\Big)\label{PP4}\\
&&+\int^t_0\!\!\!\int\big(|\p_ku_\eps|^2(\p_i(b^i_\eps-\sigma^{il}_\eps h^l_\eps)+2c_\eps)
+2\p_ku_\eps \p_kf_\eps\big)\no\\
&&+2\int^t_0\!\!\!\int\p_ku_\eps ([\p_k,\sL_\eps](u_\eps)
-\p_i(\sigma^{il}_\eps([\p_k,\sM^l](u_\eps)+\p_kg^l_\eps)))\no\\
&&+\int^t_0\!\!\!\int(h^l_\eps\p_ku_\eps+[\p_k,\sM^l_\eps](u_\eps)+\p_kg^l_\eps)^2\no\\
&&+\int^t_0\left(\int\big((\p_ku_\eps)^2\p_i\sigma^{il}_\eps+2\p_ku_\eps
([\p_k,\sM^l_\eps](u_\eps)+\p_kg^l_\eps)\big)\right)\dif B^l_t.\no
\ee
We only need to treat the trouble terms
$$
\int\p_ku_\eps\p_i(\p_ka^{ij}_\eps\p_j u_\eps)\ \mbox{ and }\
\int\p_k u_\eps\p_i(\sigma^{il}_\eps\p_k\sigma^{jl}_\eps\p_j u_\eps).
$$
The first one can be dealt with as follows:
\be
&&\int\p_ku_\eps\p_i(\p_ka^{ij}_\eps\p_j u_\eps)=-\int\p_k\p_iu_\eps\p_ka^{ij}_\eps\p_j u_\eps\no\\
&&\qquad=\int\p_iu_\eps\p_k\p_ka^{ij}_\eps\p_j u_\eps
+\int\p_iu_\eps\p_ka^{ij}_\eps\p_k\p_j u_\eps.\label{Es8}
\ee
Thus, by the symmetry of $a^{ij}_\eps$ and (\ref{P1}), we have
\be
\int\p_ku_\eps\p_i(\p_ka^{ij}_\eps\p_j u_\eps)
=\frac{1}{2}\int\p_iu_\eps\p_k\p_ka^{ij}_\eps\p_j u_\eps
\leq \ell_s\int|\p_i u_\eps|^2,\label{PP3}
\ee
where we have used that $\p_k\p_ka^{ij}_{s,\eps}\xi_i\xi_j\leq \ell_s|\xi|^2$ by (\ref{P1}).

For the second one, noticing that
$$
\p_i(\sigma^{il}_\eps\p_k\sigma^{jl}_\eps\p_j u_\eps)=\p_i\sigma^{il}_\eps\p_k\sigma^{jl}_\eps\p_j u_\eps
+\sigma^{il}_\eps\p_i\p_k\sigma^{jl}_\eps\p_j u_\eps+\sigma^{il}_\eps\p_k\sigma^{jl}_\eps\p_i\p_j u_\eps
$$
and
$$
\sigma^{il}_\eps\p_k\sigma^{jl}_\eps\p_i\p_j u_\eps=
\frac{1}{2}\p_k(\sigma^{il}_\eps\sigma^{jl}_\eps)\p_i\p_ju_\eps,
$$
as in (\ref{Es8})  and by (\ref{P3}), we have
\be
\int\p_k u_\eps\p_i(\sigma^{il}_\eps\p_k\sigma^{jl}_\eps\p_j u_\eps)&\leq& \ell_s\int|\p_i u_\eps|^2
+\frac{1}{2}\int\p_k u_\eps\p_k(\sigma^{il}_\eps\sigma^{jl}_\eps)\p_i\p_ju_\eps\no\\
&\leq& \ell_s\int|\p_i u_\eps|^2
+\frac{1}{4}\int\p_i u_\eps\p^2_k(\sigma^{il}_\eps\sigma^{jl}_\eps)\p_ju_\eps\no\\
&\leq& \ell_s\int|\p_i u_\eps|^2.\label{PP2}
\ee
By (\ref{PP4}), (\ref{PP3}), (\ref{PP2}), (\ref{P2}) and (\ref{P3}), we find that
\ce
&&\int |\p_k u_{\eps,t}|^2\leq \int |\p_k u_{\eps,0}|^2+C+\int^t_0\ell_s\int|\p_i u_\eps|^2\\
&&\qquad+\int^t_0\left(\int\big((\p_ku_\eps)^2\p_i\sigma^{il}_\eps+2\p_ku_\eps
([\p_k,\sM^l_\eps](u_\eps)+\p_kg^l_\eps)\big)\right)\dif B^l_t.
\de
Using the same method as proving (\ref{Es5}),
we may prove the following uniform estimate:
$$
\mE\left(\sup_{t\in[0,T]}\int|\p_i u_{\eps,t}|^2\right)\leq C,
$$
which then produces (\ref{PL4}).
\end{proof}

In Proposition \ref{Pro2}, certain conditions on second order derivatives of $a$ and $b$
are required. Below, we follow the idea of LeBris and Lions \cite{Le-Li} to consider a special
degenerate case so that we can weaken the assumptions on $a$ and $b$ (see (\ref{CC1}) below).
But, we need a stronger assumption than the parabolic condition (see (\ref{CC2}) below).

\bp\label{Pro3}
Let $a^{ij}$ be given as follows
\be
a^{ij}_t(\omega,x)=\hat\sigma^{il}_t(\omega,x)\hat\sigma^{jl}_t(\omega,x)\label{CC1}
\ee
such that for some $\a>1/2$
\be
|\hat\sigma^{il}_t\xi_i|^2\geq \a |\sigma^{il}_t\xi_i|^2,\ \ \forall \xi\in\mR^d.\label{CC2}
\ee
Assume also that the following conditions hold: for some $q>1$
\be
&\hat\sigma^{i\cdot},\p_i\hat\sigma^{i\cdot}\in L^2(0,T; L^\infty(\Omega; L^4_{loc}(\mR^d;l^2))),\\
&\div b,\ c\in L^1(0,T; L^\infty(\Omega\times\mR^d)),\label{CC3}\\
&\sigma^{i\cdot}, \p_i\sigma^{i\cdot}, h\in L^{2q}(0,T; L^\infty(\Omega\times\mR^d;l^2)).\label{CC4}
\ee
Then for any $u_0\in L^2(\Omega,\cF_0; L^2(\mR^d))$ and
$$
f\in L^2([0,T]\times\Omega\times\mR^d),
\ \ g\in L^2([0,T]\times\Omega\times\mR^d; l^2),
$$
there exists a generalized solution $u$ of SPDE (\ref{SPDE}) such that
\be
&&\mE\left(\mathrm{ess.}\sup_{s\in[0,T]}\int|u_{s}|^{2}\right)
+\mE\left(\int^T_0\!\!\!\int(|\hat\sigma^{il}\p_i u|^2+|\sigma^{il}\p_i u|^2)\right)\no\\
&&\qquad\qquad\leq C\left(\mE\int|u_0|^2+\mE\int^T_0\!\!\!\int|f|^2
+\mE\int^T_0\!\!\!\int|g|^2\right),\label{Ep1}
\ee
where $C$ is independent of $u_0, f$ and $g$.
\ep
\begin{proof}
Let $a^{ij}_\eps$ and $\sigma^{il}_\eps$ be defined by (\ref{de2}).
Let $\hat\sigma^{il}_\eps:=(\hat\sigma^{il}*\rho_\eps)\chi_\eps$.
As  (\ref{Es9}) and (\ref{PL2}), we have for all $\xi\in\mR^d$
$$
a^{ij}_\eps\xi_i\xi_j\geq|\hat\sigma^{il}_\eps\xi_i|^2,
\ \ a^{ij}_\eps\xi_i\xi_j\geq \a |\sigma^{il}_\eps\xi_i|^2,
$$
which implies that
\be
\sA_{a_\eps,\sigma_\eps}(\xi)=2a^{ij}_\eps\xi_i\xi_j-|\sigma^{il}_\eps\xi_i|^2
\geq \frac{2\a-1}{1+\a}
(a^{ij}_\eps\xi_i\xi_j+|\sigma^{il}_\eps\xi_i|^2)\geq \frac{2\a-1}{1+\a}(|\hat\sigma^{il}_\eps\xi_i|^2+|\sigma^{il}_\eps\xi_i|^2).\label{PL3}
\ee
In (\ref{Es1}), using the left hand side of (\ref{L6}), by (\ref{PL3}), (\ref{CC3}), (\ref{CC4})
and Young's inequality, we have
\be
\int|u_{\eps,t}|^2&\leq&\int|u_{\eps,0}|^2-\frac{2\a-1}{1+\a}\int^t_0\!\!\!\int\Big(
|\hat\sigma^{il}_\eps\p_i u_\eps|^2+|\sigma^{il}_\eps\p_i u_\eps|^2\Big)\no\\
&&+\int^t_0\left(\ell_s\int|u_\eps|^2+C\int|f_\eps|^2+C\int|g_\eps|^2\right)\dif s\no\\
&&+\int^t_0\left(\int\Big(|u_\eps|^2\p_i\sigma^{il}_\eps+2u_\eps
 (h^l_\eps u_\eps+g^l_\eps)\Big)\right)\dif B^l_s.\label{Op1}
\ee
Thus, as in proving (\ref{Es5}), we can prove that
\ce
&&\mE\left(\sup_{s\in[0,T]}\int|u_{\eps,s}|^{2}\right)
+\mE\left(\int^T_0\!\!\!\int(|\hat\sigma^{il}_\eps\p_i u_\eps|^2
+|\sigma^{il}_\eps\p_i u_\eps|^2)\right)\\
&&\qquad\qquad\leq C\left(\mE\int|u_0|^2+\mE\int^T_0\!\!\!\int|f|^2
+\mE\int^T_0\!\!\!\int|g|^2\right),
\de
where $C$ is independent of $\eps, u_0,f,g$. The existence of generalized solution now follows
by using weakly convergence method as in the proof of Proposition \ref{Th1}.
Estimate (\ref{Ep1}) now follows as in proving (\ref{Es4}) and (\ref{Es6}).
\end{proof}
\br
If $\sigma^{i\cdot}, \p_i\sigma^{i\cdot}\in
L^{2q}(0,T; L^\infty(\Omega\times\mR^d;l^2))$ are replaced by
$$
\sigma^{i\cdot}, \p_i\sigma^{i\cdot}\in
L^{2q}(0,T; L^2(\Omega;L^2_{loc}(\mR^d;l^2))),
$$
then we still have the existence of generalized
solutions. In fact, we just need to take expectations for (\ref{Op1}). Thus, we only have
\ce
&&\sup_{t\in[0,T]}\mE \int|u_{t}|^{2}
+\mE\left(\int^T_0\!\!\!\int(|\hat\sigma^{il}\p_i u|^2+|\sigma^{il}\p_i u|^2)\right)\\
&&\qquad\leq C\left(\mE\int|u_0|^2+\mE\int^T_0\!\!\!\int|f|^2
+\mE\int^T_0\!\!\!\int|g|^2\right).
\de
\er

\section{Maximal Principle and Uniqueness for SPDE}

In this section, we  prove a maximal principle for SPDEs,
which automatically produces the uniqueness of generalized solutions.

Consider the following SPDE:
\be
\dif u=(\sL u+f)\dif t+\sM^l u\dif B^l_t,\ u_0=\varphi\label{SPDE0}.
\ee
Let $u\in L^\infty(0,T;L^2(\Omega\times\mR^d))$ be a
generalized solution of (\ref{SPDE0}) in the sense of Definition \ref{Def1}.
We first make convolutions for (\ref{SPDE0}) with $\rho_\eps$ and obtain
$$
\dif (\rho_\eps*u)=[\rho_\eps*(\sL u+f)]\dif t+[\rho_\eps*(\sM^l u)]\dif B^l_t.
$$
Set
$$
u_\eps:=\rho_\eps*u,\ \ f_\eps:=\rho_\eps*f.
$$
Let $\beta\in C^2(\mR)$ be a convex function with
\be
\beta'(r),r\beta'(r)-\beta(r),\beta''(r), r^2\beta''(r)\mbox{ are bounded}.\label{Op2}
\ee
By It\^o's formula,  we have
\ce
\dif\beta(u_\eps)&=&\beta'(u_\eps) (\rho_\eps*(\sL u)+f_\eps)\dif t
+\beta'(u_\eps) (\rho_\eps*(\sM^l u))\dif B^l_t\\
&&+\frac{1}{2}\beta''(u_\eps) |\rho_\eps*(\sM^l u)|^2\dif t.
\de
Multiplying both sides by a non-negative smooth function
$\phi\in C^\infty_0(\mR^d)$  and integrating over $\mR^d$, we get
\be
\dif\int\beta(u_\eps)\phi\=
\left(\int\beta'(u_\eps)\phi(\rho_\eps*(\sL u)+f_\eps)\right)\dif t\no\\
&&+\frac{1}{2}\left(\int\beta''(u_\eps)\phi|\rho_\eps*(\sM^l u)|^2\right)\dif t\no\\
&&+\left(\int\beta'(u_\eps)\phi \rho_\eps*(\sM^l u)\right)\dif B^l_t\no\\
\=\left(\int\beta'(u_\eps)\phi(\sL u_\eps+[\rho_\eps,\sL](u)+f_\eps)\right)\dif t\no\\
&&+\frac{1}{2}\left(\int\beta''(u_\eps)\phi
|\sM^l u_\eps+[\rho_\eps,\sM^l](u)|^2\right)\dif t\no\\
&&+\left(\int\beta'(u_\eps)\phi (\sM^l u_\eps
+[\rho_\eps,\sM^l](u))\right)\dif B^l_t,\label{Cp4}
\ee
where we have used the following notation: for a differential operator $\sD$,
$$
[\rho_\eps,\sD](u):=\rho_\eps*(\sD u)-\sD(\rho_\eps*u).
$$

\br
The following two commutation relations can be verified immediately and will be used below:
for real functions $a,b,u$,
\be
&\p[\rho_\eps,a](u)=[\rho_\eps,\p a](u)+[\rho_\eps,a\p](u),\label{For1}\\
&[\rho_\eps,ab](u)=a[\rho_\eps,b\p](u)+[\rho_\eps,a](b \p u).\label{For2}
\ee
\er
Integrating both sides of (\ref{Cp4})
in time from $0$ to $t$ and using the integration by parts formula,
as in (\ref{Es1}) we further have
\be
\int\beta(u_{t,\eps})\phi
=\int\beta(u_{0,\eps})\phi+\sum_{i=1}^8J^\eps_i(t),\label{Me1}
\ee
where
\ce
J^\eps_1(t)&:=&\int^t_0\!\!\!\int\Big(-\beta''(u_\eps)\phi
\frac{\sA_{a,\sigma}(\nabla u_\eps)}{2}
+\beta'(u_\eps)\phi f_\eps\Big),\no\\
J^\eps_2(t)&:=&\int^t_0\!\!\!\int\beta'(u_\eps)\phi[\rho_\eps,\sL](u),\no\\
J^\eps_3(t)&:=&\int^t_0\!\!\!\int\beta(u_\eps)[\p_j(\p_i\phi a^{ij})
-\p_i\phi b^i+c\phi],\no\\
J^\eps_4(t)&:=&\int^t_0\!\!\!\int(u_\eps\beta'(u_\eps)-\beta(u_\eps))
[\phi\p_ib^i+c\phi],\no\\
J^\eps_5(t)&:=&\int^t_0\!\!\!\int\beta''(u_\eps)\phi\sigma^{il}\p_iu_\eps
[\rho_\eps,\sM^l](u),\no\\
J^\eps_6(t)&:=&\int^t_0\!\!\!\int\beta''(u_\eps)\phi\sigma^{il}\p_iu_\eps
 h^l u_\eps,\no\\
J^\eps_7(t)&:=&\frac{1}{2}\int^t_0\!\!\!\int\beta''(u_\eps)\phi
(h^l  u_\eps+[\rho_\eps,\sM^l](u))^2,\no\\
J^\eps_8(t)&:=&\int^t_0\left(\int\beta'(u_\eps)\phi (\sM^l u_\eps
+[\rho_\eps,\sM^l](u))\right)\dif B^l_s.
\de

We want to take limits $\eps\downarrow 0$.
For this aim, we need the following key commutation lemma of DiPerna-Lions \cite{Di-Li}.
For the reader's convenience, a detailed proof is provided in the appendix.
\bl\label{Le2}
For $j=1,2,3$, let $p_j\in[1,\infty]$ and $q_j\geq\frac{p_j}{p_j-1}$. We are given
$$
u\in L^{p_1}(0,T; L^{p_2}(\Omega; L^{p_3}_{loc}(\mR^d))),
\ \ c\in L^{q_1}(0,T;L^{q_2}(\Omega;L^{q_3}_{loc}(\mR^d)))
$$
and for $i=1,\cdots,d$,
$$
b^i\in L^{q_1}(0,T;L^{q_2}(\Omega;W^{1,q_3}_{loc}(\mR^d))).
$$
Let $r_j\in[1,\infty)$ be given by $\frac{1}{r_j}=\frac{1}{p_j}+\frac{1}{q_j}, j=1,2,3$. Then,
\be
[\rho_\eps,b^i\p_i](u)\stackrel{\eps\to 0}{\longrightarrow}0
\ \mbox{ in $L^{r_1}(0,T;L^{r_2}(\Omega; L^{r_3}_{loc}(\mR^d)))$}\label{Lm1}
\ee
and
\be
[\rho_\eps,c](u)\stackrel{\eps\to 0}{\longrightarrow}0
\ \mbox{ in $L^{r_1}(0,T;L^{r_2}(\Omega, L^{r_3}_{loc}(\mR^d)))$}.\label{Lm2}
\ee
Moreover, if
$$
u\in L^{p_1}(0,T; L^{p_2}(\Omega; W^{1,p_3}_{loc}(\mR^d))),\ \
b^i\in L^{q_1}(0,T;L^{q_2}(\Omega;L^{q_3}_{loc}(\mR^d))),
$$
then (\ref{Lm1}) still holds.
\el

We first treat the terms $J^\eps_2,J^\eps_3,J^\eps_4$.
\bl\label{Le4}
Let $u\in L^2([0,T]\times\Omega; W^{1,2}_{loc}(\mR^d))\cap
L^\infty(0,T; L^2(\Omega; L^2_{loc}(\mR^d)))$
and assume that
\be
a^{ij}\in L^\infty([0,T]\times\Omega;L^\infty_{loc}(\mR^d))
\cup L^2(0,T; L^\infty(\Omega;W^{1,\infty}_{loc}(\mR^d)))\label{L11}
\ee
and
\be
b^i\in L^2(0,T; L^2(\Omega; L^2_{loc}(\mR^d))),\ \div b, c
\in L^1(0,T; L^2(\Omega; L^2_{loc}(\mR^d)))\label{L33}
\ee
or
\be
b^i\in L^1(0,T; L^2(\Omega; W^{1,2}_{loc}(\mR^d))),\ c
\in L^1(0,T; L^2(\Omega; L^2_{loc}(\mR^d)))\label{L333}
\ee
hold. Then, we have
\be
\lim_{\eps\to 0}\mE|J^\eps_2(t)|=0\label{Lm11}
\ee
and in $L^1(\Omega)$
\be
J^\eps_3(t)&\stackrel{\eps\to 0}{\longrightarrow}&\int^t_0\!\!\!
\int\beta(u)[\p_j(\p_i\phi a^{ij})
-\p_i\phi b^i+c\phi],\label{Lp1}\\
J^\eps_4(t)&\stackrel{\eps\to 0}{\longrightarrow}&\int^t_0\!\!\!\int(u\beta'(u)-\beta(u))
[\phi\p_ib^i+c\phi],\label{Lp2}
\ee
where for (\ref{Lp1}), we also need the assumption $\p_ja^{ij}\in L^1(0,T; L^2(\Omega;L^2_{loc}(\mR^d)))$.
\el
\begin{proof}
By (\ref{For1}), we have
$$
[\rho_\eps,\sL](u)=\p_i[\rho_\eps,a^{ij}\p_j](u)+[\rho_\eps,\p_ib^i](u)
+[\rho_\eps,b^i\p_i](u)+[\rho_\eps,c](u).
$$
Thus, we may write
\be
J_2^\eps(t)
&=&-\int^t_0\!\!\!\int\beta''(u_\eps)\p_iu_\eps\phi[\rho_\eps,a^{ij}\p_j](u)\no\\
&&-\int^t_0\!\!\!\int\beta'(u_\eps)\p_i\phi[\rho_\eps,a^{ij}\p_j](u)
+\int^t_0\!\!\!\int\beta'(u_\eps)\phi[\rho_\eps,\p_ib^i](u)\no\\
&&+\int^t_0\!\!\!\int\beta'(u_\eps)\phi[\rho_\eps,b^i\p_i](u)
+\int^t_0\!\!\!\int\beta'(u_\eps)\phi[\rho_\eps,c](u)\no\\
&=:&J^\eps_{21}(t)+J^\eps_{22}(t)+J^\eps_{23}(t)+J^\eps_{24}(t)+J^\eps_{25}(t).\label{Lp5}
\ee
Let $Q:=\mathrm{supp}(\phi)$. By H\"older's inequality, we have
\be
\mE|J^\eps_{21}(t)|\leq C\left(\mE\int^t_0\!\!\!\int_Q|\p_iu_\eps|^2\right)^{1/2}
\left(\mE\int^t_0\!\!\!\int_Q|[\rho_\eps,a^{ij}\p_j](u)|^2\right)^{1/2},\label{L5}
\ee
which converges to zero as $\eps\to 0$ by (\ref{L11})  and (\ref{Lm1}) or
the second conclusion of Lemma \ref{Le2}.
Similarly,
$$
\mE|J^\eps_{22}(t)|\stackrel{\eps\to 0}{\longrightarrow} 0.
$$
Moreover, by Lemma \ref{Le2}, (\ref{Op2}) and (\ref{L33}) or (\ref{L333}), we also have
\be
\mE|J^\eps_{23}(t)|+\mE|J^\eps_{24}(t)|+\mE|J^\eps_{25}(t)|
\stackrel{\eps\to 0}{\longrightarrow} 0.\label{PP5}
\ee
Limit (\ref{Lm11}) now follows. Limits (\ref{Lp1}) and (\ref{Lp2}) are easy by (\ref{Op2})
and the dominated convergence theorem.
\end{proof}

Next, we look at the term $J^\eps_5$.
\bl\label{Le5}
Let $u\in L^2([0,T]\times\Omega; W^{1,2}_{loc}(\mR^d))\cap
L^\infty(0,T; L^2(\Omega; L^2_{loc}(\mR^d)))$
and assume that
\be
\sigma^{i\cdot}\in L^\infty([0,T]\times\Omega;L^\infty_{loc}(\mR^d;l^2))), \ \
 h\in L^2(0,T; L^\infty(\Omega;L^\infty_{loc}(\mR^d;l^2)))\label{L4}
\ee
hold. Then, we have
\be
\lim_{\eps\to 0}\mE |J^\eps_5(t)|=0.\label{Lm111}
\ee
\el
\begin{proof}
In view of
$$
[\rho_\eps,\sM^l](u)=[\rho_\eps,\sigma^{il}\p_i](u)+[\rho_\eps,h^l](u),
$$
by (\ref{Lm1}), (\ref{Lm2}) and (\ref{L4}), one sees that
\be
[\rho_\eps,\sM^l](u)\stackrel{\eps\to 0}{\longrightarrow} 0
\mbox{ in $L^2(0,T; L^2(\Omega;L^2_{loc}(\mR^d)))$}.\label{PP6}
\ee
Limit (\ref{Lm111}) now follows by (\ref{L4}) and $u\in L^2([0,T]\times\Omega; W^{1,2}_{loc}(\mR^d))$.
\end{proof}
\br
In Lemmas \ref{Le4} and \ref{Le5}, if we assume
$$
\p_i u\in L^\infty(0,T;L^2(\Omega; L^2_{loc}(\mR^d))),
$$
then  the conditions on $a$ and $\sigma$ in (\ref{L11}) and (\ref{L4})  can be replaced by
$$
a^{ij}\in L^1(0,T; L^\infty(\Omega;L^\infty_{loc}(\mR^d))),\ \
\sigma^{i\cdot}\in L^2(0,T; L^\infty(\Omega;L^\infty_{loc}(\mR^d;l^2))).
$$
\er
We first prove:
\bp\label{Pro4}
Assume that (\ref{L11}), (\ref{L33}), (\ref{L4}) and the following conditions hold:
\be
&c^+ \in L^1(0,T; L^\infty(\Omega\times\mR^d)),\label{L3}\\
&\p_i\sigma^{i\cdot},\ \p_k h\in L^2(0,T; L^\infty(\Omega;L^2_{loc}(\mR^d;l^2))).\label{L44}
\ee
Let $u\in L^\infty(0,T;L^2(\Omega\times\mR^d))$ be a generalized solution of (\ref{SPDE0}) satisfying
\be
\p_i u\in L^2([0,T]\times\Omega; L^2_{loc}(\mR^d)).\label{L2}
\ee
\begin{enumerate}[{\rm (I)}]
\item
If $f\geq 0$ and $u_0\geq 0$ and one of the following conditions holds
\be
\frac{|a^{ij}|}{1+|x|^2},\ \frac{|\p_j a^{ij}|}{1+|x|},\ \frac{|b^i|}{1+|x|}
\in L^1(0,T; L^2(\Omega\times\mR^d)),\label{C3}
\ee
\be
\frac{|a^{ij}|}{1+|x|^2},\ \frac{|\p_j a^{ij}|}{1+|x|},
\ \frac{|b^i|}{1+|x|}\in L^1(0,T; L^\infty(\Omega\times\mR^d)),\label{C33}
\ee
then for $(\dif t\times P\times\dif x)$-almost
all $(t,\omega,x)\in[0,T]\times\Omega\times\mR^d$
$$
u_t(\omega,x)\geq 0.
$$
\item If  $u_0\in L^1(\Omega,\cF_0;L^1(\mR^d))$,
$f\in L^1([0,T]\times\Omega\times\mR^d)$ and (\ref{C3}) together with the following condition holds:
\be
\frac{\|\sigma^{i\cdot}\|_{l^2}}{1+|x|}\in L^2(0,T;L^\infty(\Omega;L^2(\mR^d))),
\ \ \p_i\sigma^{i\cdot}, h\in L^2(0,T; L^\infty(\Omega\times\mR^d;l^2)),\label{C30}
\ee
then
\be
\mE\left(\mathrm{ess.}\sup_{t\in[0,T]}\int|u_t|\right)
\leq C\mE\int|u_0|+C\mE\int^T_0\!\!\!\int |f|,\label{Es7}
\ee
where the constant $C$ only depends on
$\|\p_i\sigma^{i\cdot}\|_{L^2(0,T; L^\infty(\Omega\times\mR^d;l^2))}$,
$\|h\|_{L^2(0,T; L^\infty(\Omega\times\mR^d;l^2))}$ and\\
$\|c^+\|_{L^1(0,T; L^\infty(\Omega\times\mR^d))}.$
\end{enumerate}
\ep
\begin{proof}
Using the integration by parts formula
and by (\ref{L4}), (\ref{L44}) and the dominated convergence theorem, we have
\be
J^\eps_6(t)&=&\int^t_0\!\!\!\int(\beta(u_\eps)-u_\eps\beta'(u_\eps))\p_i(\phi\sigma^{il}h^l)\no\\
&\stackrel{\eps\to 0}{\longrightarrow}&\int^t_0\!\!\!
\int(\beta(u)-u\beta'(u))\p_i(\phi\sigma^{il}h^l)\  \mbox{ in $L^1(\Omega)$},\label{PP9}
\ee
and by (\ref{PP6}),
\be
J^\eps_7(t)\stackrel{\eps\to 0}{\longrightarrow}\int^t_0\!\!\!
\int\beta''(u)\phi(h^l  u)^2\  \mbox{ in $L^1(\Omega)$.}\label{PP7}
\ee
Moreover, we also have
\ce
J^\eps_8(t)&=&\int^t_0\left(\int\Big(\beta(u_\eps)\p_i(\sigma^{il}\phi)
+\beta'(u_\eps) \phi (h^l u_\eps +[\rho_\eps,\sM^l](u))\Big)\right)\dif B^l_s\\
&\stackrel{\eps\to 0}{\longrightarrow}&\int^t_0
\left(\int\Big(\beta(u)\p_i(\sigma^{il}\phi)+\beta'(u)\phi h^l u\Big)\right)\dif B^l_s\
\mbox{ in $L^2(\Omega)$,}
\de
where the above stochastic integral is a continuous $L^2$-martingale.

Now taking limits $\eps\to 0$ for (\ref{Me1}) and summarizing the above limits,
we arrive at
\be
\int\beta(u_t)\phi&\leq&\int\beta(u_0)\phi+\int^t_0\!\!\!\int(u\beta'(u)-\beta(u))
(\phi\p_ib^i-\p_i(\phi\sigma^{il} h^l)+c\phi)\no\\
&&+\int^t_0\!\!\!\int\beta(u)(\p_j(\p_i\phi a^{ij})-\p_i\phi b^i+c\phi)\no\\
&&+\int^t_0\!\!\!\int\beta'(u)\phi f+\frac{1}{2}\int^t_0\!\!\!\int\beta''(u)\phi(h^l  u)^2\no\\
&&+\int^t_0\left(\int\Big(\beta(u)\p_i(\sigma^{il}\phi)
+\beta'(u)\phi h^l u)\right)\dif B^l_s.\label{PP10}
\ee
(I). Let $\beta(r)=\beta_\delta(r)=\frac{\sqrt{r^2+\delta}-r}{2}$ in (\ref{PP10}).
By simple calculations, we have
$$
\lim_{\delta\downarrow 0}\beta_\delta(r)=-(0\wedge r):=r^{-},\ \ \beta'_\d(r)\leq 0,
$$
and
$$
|r\beta_\delta'(r)-\beta_\delta(r)|\leq\frac{\sqrt{\delta}}{2},\ \
|r^2\beta''_\delta(r)|\leq\frac{\sqrt{\delta}}{2}.
$$
Taking expectations for  (\ref{PP10}) and letting $\delta\to 0$, by (\ref{L4}),
(\ref{L3}), (\ref{L44}) and $f\geq 0$, $u_0\geq 0$, we get
$$
\mE\int u^-_t\phi\leq\mE\int^t_0\!\!\!\int
(\p_j(\p_i\phi a^{ij})-\p_i\phi b^i) u^-_s+\int^t_0\ell_s\mE\int\phi u^-_s,
$$
which yields by Gronwall's inequality,
\be
\mE\int u^-_t\phi\leq C\mE\int^t_0\!\!\!\int (\p_j(\p_i\phi a^{ij})-\p_i\phi b^i) u^-_s.
\label{L7}
\ee

{\bf Case (\ref{C3}):}
Let $\chi_n(x)=\chi(x/n)$ be a cutoff function with the same $\chi$ as in Section 3.
We choose in (\ref{L7})
$$
\phi(x)=\chi_n(x).
$$
Noting that
$$
|\p_i\chi_n(x)|\leq\frac{C1_{n\leq|x|\leq2n}}{n}\leq\frac{C1_{|x|\geq n}}{1+|x|}
$$
and
$$
|\p_i\p_j\chi_n(x)|\leq\frac{C1_{n\leq|x|\leq2n}}{n^2}\leq\frac{C1_{|x|\geq n}}{1+|x|^2},
$$
we have
\be
\mE\int u^-_t\chi_n&\leq& C\mE\int^t_0\!\!\!\int_{|x|\geq n}
\left(\frac{|a^{\cdot\cdot}|}{1+|x|^2}+
\frac{|\p_j a^{\cdot j}|}{1+|x|}+\frac{|b|}{1+|x|}\right) u^-_s\no\\
&\leq& C\int^t_0\left(\mE\int_{|x|\geq n} \left(\frac{|a^{\cdot\cdot}|}{1+|x|^2}+
\frac{|\p_j a^{\cdot j}|}{1+|x|}+\frac{|b|}{1+|x|}\right)^2\right)^{1/2}.\label{Lp4}
\ee
Letting $n\to\infty$ and by Fatou's lemma and (\ref{C3}), we obtain
$$
\mE\int u^-_t=0.
$$

{\bf Case (\ref{C33}):}
Let $\lambda(x):=(1+|x|^2)^{-d}$ be a weight function and choose in (\ref{L7})
$$
\phi(x)=\phi_n(x)=\lambda(x)\chi_n(x).
$$
Noting that
$$
|\p_i\phi_n(x)|\leq\lambda(x)\frac{C 1_{|x|\geq n}}{1+|x|}+
\frac{C\phi_n(x)}{1+|x|}
$$
and
$$
|\p_i\p_j\phi_n(x)|\leq\lambda(x)\frac{C 1_{|x|\geq n}}{1+|x|^2}
+\frac{C\phi_n(x)}{1+|x|^2},
$$
we have
\ce
\mE\int u^-_t\phi_n&\leq& C\int^t_0\mE\int_{|x|\geq n} \left(\frac{|a^{\cdot\cdot}|}{1+|x|^2}+
\frac{|\p_j a^{\cdot j}|}{1+|x|}+\frac{|b|}{1+|x|}\right) u^-_s\lambda\\
&& +C\int^t_0\mE\int \left(\frac{|a^{\cdot\cdot}|}{1+|x|^2}+
\frac{|\p_j a^{\cdot j}|}{1+|x|}+\frac{|b|}{1+|x|}\right) u^-_s\phi_n\\
&\leq& C\int^t_0\ell_s\mE\int_{|x|\geq n}  u^-_s\lambda
+\int^t_0\ell_s\mE\int u^-_s\phi_n.
\de
By Gronwall's inequality and letting $n\to\infty$, we get
$$
\mE\int u^-_t\lambda\leq C\lim_{n\to\infty}
\int^t_0\ell_s\mE\int_{|x|\geq n}  u^-_s\lambda=0.
$$
(II). Let $\beta_\delta(r)=\sqrt{r^2+\delta}$ in (\ref{PP10}).
By elementary calculations, we know
$$
\lim_{\delta\downarrow 0}\beta_\delta(r)=|r|, \ \ |\beta'_\delta(r)|\leq 1,
$$
and
$$
|r\beta_\delta'(r)-\beta_\delta(r)|\leq\sqrt{\delta},\ \
|r^2\beta''_\delta(r)|\leq\sqrt{\delta}.
$$
Letting $\delta\to 0$, as above we find that
\be
\int|u_t|\phi&\leq&\int|u_0|\phi
+\int^t_0\!\!\!\int|u|(\p_j(\p_i\phi a^{ij})-\p_i\phi b^i+c\phi)\no\\
&&+\int^t_0\!\!\!\int\phi |f|
+\int^t_0\left(\int|u|\big(\p_i(\sigma^{il}\phi)+\phi h^l\big)\right)\dif B^l_s,\label{PP11}
\ee
where we have used that
$$
\int^t_0\left(\int\big(\beta_\delta(u)\p_i(\sigma^{il}\phi)
+u\beta'_\delta(u)h^l\phi\big)\right)\dif B^l_s\to
\int^t_0\left(\int|u|\big(\p_i(\sigma^{il}\phi)+h^l\phi\big)\right)\dif B^l_s
$$
in $L^2(\Omega)$ as $\delta\to 0$.

Set
$$
\Phi_t:=\mathrm{ess.}\sup_{s\in[0,t]}\int|u_s|\phi.
$$
Now taking the essential supremum for both sides of (\ref{PP11})
in time $t$ and by Burkholder's inequality,   we have
\ce
\mE\Phi_t&\leq&\mE\int|u_0|\phi
+\mE\int^t_0\!\!\!\int|u|(|\p_j(\p_i\phi a^{ij})-\p_i\phi b^i|+c^+\phi)\\
&&+\mE\int^t_0\!\!\!\int\phi |f|
+C\mE\left(\int^t_0\left|\int|u|(|\p_i(\sigma^{il}\phi)|+\phi|h^l|)\right|^2\dif s\right)^{1/2}.
\de
The last term denoted by $\sI$ is controlled as follows: by (\ref{L4}), (\ref{L44})
and Young's inequality,
\ce
\sI&\leq& C\mE\left(\int^t_0\left|\int|u||\sigma^{il}\p_i\phi|\right|^2\dif s\right)^{1/2}
+\mE\left(\int^t_0\ell_s\left(\int|u_s|\phi\right)^2\dif s\right)^{1/2}\\
&\leq&C\mE\left(\int^t_0\left|\int|u||\sigma^{il}\p_i\phi|\right|^2\dif s\right)^{1/2}
+\frac{1}{2}\mE\Phi_t+\int^t_0\ell_s\mE\Phi_s\dif s.
\de
Thus, we get
\ce
\mE\Phi_t&\leq&\mE\int|u_0|\phi
+\mE\int^T_0\!\!\!\int|u|(|\p_j(\p_i\phi a^{ij})|+|\p_i\phi b^i|)
+\int^t_0\ell_s\mE\Phi_s\\
&&+\mE\int^T_0\!\!\!\int  |f|\phi+\frac{1}{2}\mE\Phi_t+
C\mE\left(\int^T_0\left|\int|u||\sigma^{il}\p_i\phi|\right|^2\dif s\right)^{1/2},
\de
which yields by Gronwall's inequality,
\be
\mE\Phi_T&\leq& C\mE\int|u_0|\phi
+C\mE\int^T_0\!\!\!\int|u|(|\p_j(\p_i\phi a^{ij})|+|\p_i\phi b^i|)\no\\
&&+C\mE\int^T_0\!\!\!\int |f|\phi
+C\mE\left(\int^t_0\left|\int|u||\sigma^{il}\p_i\phi|\right|^2\dif s\right)^{1/2}.
\label{Cp2}
\ee
Choosing $\phi=\chi_n$ and letting $n\to\infty$, as (\ref{Lp4}), we get by (\ref{C3}), (\ref{C30})
 and the dominated convergence theorem,
$$
\mE\left(\mathrm{ess.}\sup_{s\in[0,T]}\int|u_s|\right)
\leq C\mE\int|u_0|+C\mE\int^T_0\!\!\!\int |f|.
$$
If we check the above proof, we find that the constant
$C$ only depends on the following three quantities:
$$
\|\p_i\sigma^{i\cdot}\|_{L^2(0,T; L^\infty(\Omega\times\mR^d;l^2))},
\|h\|_{L^2(0,T; L^\infty(\Omega\times\mR^d;l^2))},
\|c^+\|_{L^1(0,T; L^\infty(\Omega\times\mR^d))}.
$$
The proof is complete.
\end{proof}

In the case of  $a$ and $\sigma$ independent of $x$, we have the following simple result.
\bp\label{Pro6}
Let $a,\sigma$ be independent of $x$.
Assume that the following conditions hold:
\ce
&a^{ij}\in L^1(0,T; L^\infty(\Omega)),\ \
\sigma^{i\cdot}\in L^2(0,T;L^\infty(\Omega; l^2)),\\
&b^i\in L^1(0,T; L^2(\Omega; W^{1,2}_{loc}(\mR^d))),\\
&c\in L^1(0,T; L^2(\Omega; L^2_{loc}(\mR^d))),
c^+ \in L^1(0,T; L^\infty(\Omega\times\mR^d)),\\
&h,\p_k h\in L^2(0,T; L^\infty(\Omega;L^\infty_{loc}(\mR^d;l^2))).
\de
Let $u\in L^\infty(0,T;L^2(\Omega\times\mR^d))$ be a generalized solution of (\ref{SPDE0}).
\begin{enumerate}[{\rm (I)}]
\item
If $f\geq 0$, $u_0\geq 0$  and the following condition holds
$$
\frac{|b^i|}{1+|x|} \in L^1(0,T; L^2(\Omega\times\mR^d))
\cup L^1(0,T; L^\infty(\Omega\times\mR^d)),
$$
then for $(\dif t\times P\times\dif x)$-almost
all $(t,\omega,x)\in[0,T]\times\Omega\times\mR^d$
$$
u_t(\omega,x)\geq 0.
$$
\item If  $u_0\in L^1(\Omega,\cF_0;L^1(\mR^d))$,
$f\in L^1([0,T]\times\Omega\times\mR^d)$ and  the following condition holds
$$
\frac{|b^i|}{1+|x|} \in L^1(0,T; L^2(\Omega\times\mR^d)),
$$
then
$$
\mE\left(\mathrm{ess.}\sup_{t\in[0,T]}\int|u_t|\right)\leq
C\mE\int|u_0|+C\mE\int^T_0\!\!\!\int |f|.
$$
\end{enumerate}
\ep
\begin{proof}
Noticing that in this case
$$
[\rho_\eps,\sL](u)=[\rho_\eps,\p_ib^i](u)+[\rho_\eps,b^i\p_i](u)+[\rho_\eps,c](u)
$$
and
$$
[\rho_\eps,\sM^l](u)=[\rho_\eps,h^l](u).
$$
we can repeat the proof given in Proposition \ref{Pro4} to conclude the result. We omit the details.
\end{proof}

\bp\label{Pro5}
Let $a^{ij}$ be given as follows
$$
a^{ij}_t(\omega,x)=\hat\sigma^{il}_t(\omega,x)\hat\sigma^{jl}_t(\omega,x)
$$
such that for some $\a>1/2$
\be
|\hat\sigma^{il}\xi_i|^2\geq \a|\sigma^{il}\xi_i|^2.\label{Lm4}
\ee
Assume that the following conditions hold:
\be
&\hat\sigma^{i\cdot},\p_i\hat\sigma^{i\cdot}, \sigma^{i\cdot},
\p_i\sigma^{i\cdot}\in
L^2(0,T; L^\infty(\Omega; L^\infty_{loc}(\mR^d; l^2))),\label{LLL1}\\
&b^i\in L^1(0,T; L^2(\Omega; W^{1,2}_{loc}(\mR^d))),\label{LLL33}\\
&c\in L^1(0,T; L^2(\Omega; L^2_{loc}(\mR^d))),
c^+ \in L^1(0,T; L^\infty(\Omega\times\mR^d)),\\
& h\in L^2(0,T; L^\infty(\Omega;L^\infty_{loc}(\mR^d;l^2))).\label{LLL4}
\ee
Let $u\in L^\infty(0,T;L^2(\Omega\times\mR^d))$ be a
generalized solution of (\ref{SPDE}) satisfying
\be
\hat\sigma^{i\cdot}\p_i u, \sigma^{i\cdot}\p_i u
\in L^2([0,T]\times\Omega; L^2_{loc}(\mR^d;l^2)).\label{Lm5}
\ee
Then, the same conclusions (I) and (II) of Proposition \ref{Pro4} hold.
\ep
\begin{proof}
Following the proof of Proposition \ref{Pro4}, by (\ref{Lm4}), as (\ref{PL3}) we have
\be
\int\beta(u_{t,\eps})\phi&\leq&\int\beta(u_{0,\eps})\phi
-\frac{2\a-1}{\a+1}\int^t_0\!\!\!\int\beta''(u_\eps)\phi(|\hat\sigma^{il}\p_iu_\eps|^2+
|\sigma^{il}\p_iu_\eps|^2)\no\\
&&+\int^t_0\!\!\!\int\beta'(u_\eps)\phi f_\eps+\sum_{i=2}^8|J^\eps_i(t)|,\label{Lm33}
\ee
where $J^\eps_i(t)$ are the same as in (\ref{Me1}).
Checking the proof of Lemma \ref{Le4},
we need to give different treatments for $J^\eps_{21}$ and $J^\eps_{22}$ in (\ref{Lp5}).
By (\ref{For1}), (\ref{For2}) and Young's inequality, we have for any $\delta>0$,
\be
J^\eps_{21}(t)&=&-\int^t_0\!\!\!\int\beta''(u_\eps)\phi\hat\sigma^{il}\p_iu_\eps
[\rho_\eps,\hat\sigma^{jl}\p_j](u)\no\\
&&-\int^t_0\!\!\!\int\beta''(u_\eps)\phi \p_iu_\eps
[\rho_\eps,\hat\sigma^{il}](\hat\sigma^{jl}\p_ju)\no\\
&\leq&\delta\int^t_0\!\!\!\int\beta''(u_\eps)\phi|\hat\sigma^{il}\p_iu_\eps|^2\no\\
&&+C_\delta\int^t_0\!\!\!\int\beta''(u_\eps)\phi|[\rho_\eps,\hat\sigma^{jl}\p_j](u)|^2\no\\
&&+\int^t_0\!\!\!\int\beta'(u_\eps)\p_i\phi[\rho_\eps,\hat\sigma^{il}]
(\hat\sigma^{jl}\p_ju)\no\\
&&+\int^t_0\!\!\!\int\beta'(u_\eps)\phi[\rho_\eps,\p_i\hat\sigma^{il}]
(\hat\sigma^{jl}\p_ju)\no\\
&&+\int^t_0\!\!\!\int\beta'(u_\eps)\phi
[\rho_\eps,\hat\sigma^{il}\p_i](\hat\sigma^{jl}\p_ju).\label{Lm6}
\ee
By (\ref{Lm5}), (\ref{LLL1}) and Lemma \ref{Le2}, except for the first term, the other terms tend to
zero in $L^1(\Omega)$ as $\eps\to 0$. As for $J^\eps_{22}$ in  (\ref{Lp5}), we can treat it
in the same way as above, and have
$$
\lim_{\eps\to 0}\mE|J^\eps_{22}(t)|=0.
$$

For $J^\eps_5$ and $J^\eps_6$, by Young's inequality, we have for any $\delta>0$,
\be
J^\eps_5(t)+J^\eps_6(t)
&\leq&\d\int^t_0\!\!\!\int\beta''(u_\eps)\phi|\sigma^{il}\p_iu_\eps|^2
+C_\delta\int^t_0\!\!\!\int\beta''(u_\eps)\phi |h^l u_\eps|^2\no\\
&&+C_\delta\int^t_0\!\!\!\int\beta''(u_\eps)\phi |[\rho_\eps,\sM^l](u)|^2.\label{Lm7}
\ee
By (\ref{LLL1}), (\ref{LLL4}) and Lemma \ref{Le4}, the last term goes to zero as $\eps\to 0$.
Substituting (\ref{Lm6}) and (\ref{Lm7}) into (\ref{Lm33}), taking $\delta$ small enough and
letting $\eps\to 0$, we obtain
\ce\int\beta(u_t)\phi&\leq&\int\beta(u_0)\phi+\int^t_0\!\!\!\int(u\beta'(u)-\beta(u))
(\phi\p_ib^i+c\phi)\no\\
&&+\int^t_0\!\!\!\int\beta(u)(\p_j(\p_i\phi a^{ij})-\p_i\phi b^i+c\phi)\no\\
&&+\int^t_0\!\!\!\int\beta'(u)\phi f+C_\delta\int^t_0\!\!\!\int\beta''(u)\phi(h^l  u)^2\no\\
&&+\int^t_0\left(\int\Big(\beta(u)\p_i(\sigma^{il}\phi)
+\beta'(u)\phi h^l u)\right)\dif B^l_s.
\de
Thus, we can repeat the proof  of Proposition \ref{Pro4}. The details are omitted.
\end{proof}

\section{$L^1$-Integrability and Weak Continuity of Generalized Solutions}

Although we have already proved the $L^1$-integrability of generalized
solutions in the previous section under (\ref{C3}) and (\ref{C30}), we still hope to
get the $L^1$-integrability under (\ref{C33}). We now return to the construction of
generalized solutions and use estimate (\ref{Es7}) to prove the $L^1$-integrability
of the constructed solutions in Section 3.
Moreover, we shall also study the weak continuity of generalized
solutions in $L^2(\mR^d)$ and $L^1(\mR^d)$.

As in Section 3, we start from approximation equation (\ref{App}).
Instead of there, we use the following approximation for $b$ as used in (\ref{de2}):
$$
b^i_{t,\eps}:=(b^i_t*\rho_\eps)\chi_\eps.
$$
\bl\label{Le3}
Assume that
$$
\frac{|b|}{1+|x|}, \ \p_i b^i\in L^1(0,T;L^\infty(\Omega\times\mR^d)).
$$
Then for some $\ell_t\in L^1(0,T)$,
\be
\sup_{\eps\in(0,1)}\|\p_ib^i_{t,\eps}\|_{L^\infty(\Omega\times\mR^d)}\leq\ell_t.\label{PP1}
\ee
\el
\begin{proof}
Note that
$$
\p_ib^i_{t,\eps}=\p_i(b^i_t*\rho_\eps)\chi_\eps+(b^i_t*\rho_\eps)\p_i\chi_\eps.
$$
It is clear that
$$
\|\p_i(b^i_t*\rho_\eps)\chi_\eps\|_{L^\infty(\Omega\times\mR^d)}\leq
C\|\p_ib^i_t\|_{L^\infty(\Omega\times\mR^d)}.
$$
Observing (\ref{PL5}) and
$$
\mathrm{supp}(\rho_\eps)\subset\{x\in\mR^d: |x|\leq \eps\},
$$
we have for $\eps\in(0,1)$
\ce
|(b^i_t*\rho_\eps)(x)\p_i\chi_\eps(x)|&\leq&
C\eps 1_{[1/\eps,2/\eps]}(|x|)\int |b_t(y)|\rho_\eps(x-y)\dif y\\
&\leq&C\eps \int_{1/\eps-1\leq|y|\leq 2/\eps+1} |b_t(y)|\rho_\eps(x-y)\dif y\\
&\leq&C\sup_{y}\frac{|b_t(y)|}{1+|y|} \int_{1/\eps-1\leq|y|\leq 2/\eps+1}
\eps(1+|y|)\rho_\eps(x-y)\dif y\\
&\leq&C\sup_{y}\frac{|b_t(y)|}{1+|y|}.
\de
Hence,
$$
\|(b^i_t*\rho_\eps)\p_i\chi_\eps\|_{L^\infty(\Omega\times\mR^d)}
\leq C\left\|\frac{|b_t|}{1+|x|}\right\|_{L^\infty(\Omega\times\mR^d)}.
$$
The desired estimate follows.
\end{proof}

We have:
\bp\label{Pro1}
Keep the same assumptions as in Proposition \ref{Th1} and assume
$$
\frac{|b|}{1+|x|}\in L^1(0,T;L^\infty(\Omega\times\mR^d)).
$$
If $g^l\equiv 0$, $u_0\in L^1(\Omega,\cF_0;L^1(\mR^d))$ and
$f\in L^1([0,T]\times\Omega\times\mR^d)$, then
the generalized solution in Proposition \ref{Th1} satisfies
\be
\mE\left(\mathrm{ess.}\sup_{t\in[0,T]}\int|u_t|\right)
\leq C\mE\int |u_{0}|+C\mE\int^T_0\!\!\!\int|f|.\label{Es77}
\ee
\ep
\begin{proof}
Consider approximation equation (\ref{App}).
Since all the coefficients have supports contained in the ball of radius $1/\eps$,
all of the conditions  in Proposition \ref{Pro4} are satisfied. Thus,
by (\ref{Es7}), we have the following uniform estimate:
\be
\mE\left(\mathrm{ess.}\sup_{s\in[0,T]}\int|u_{\eps,s}|\right)\leq C\mE\int |u_{0}|
+C\mE\int^T_0\!\!\!\int|f|,\label{Es09}
\ee
where $C$ is independent of $\eps$.
Now, following the proof of Proposition \ref{Th1}, let
$\sD\subset C^\infty_0(\mR^d)$ be a countable and dense subset of $L^2(\mR^d)$.
Then
\ce
\sup_{t\in S_\omega}\|u_t(\omega)\|_{L^1}&=&\sup_{t\in S_\omega}\|\sqrt{|u_t(\omega)|}\|_{L^2}^2
=\sup_{t\in S_\omega}\left(\sup_{\phi\in\sD}\frac{1}{\|\phi\|_{L^2}}\int\sqrt{|u_t(\omega)|}\phi\right)^2\\
&=&\sup_{t\in S_\omega}\left(\sup_{\phi\in\sD}\frac{1}{\|\phi\|_{L^2}}\lim_{n\to\infty}
\int\sqrt{|\bar u_{\eps_n,t}(\omega)|}\phi\right)^2\\
&\leq&\varliminf_{n\to\infty}\left(\sup_{t\in S_\omega}\sup_{\phi\in\sD}\frac{1}{\|\phi\|_{L^2}}
\int\sqrt{|\bar u_{\eps_n,t}(\omega)|}\phi\right)^2\\
&=&\varliminf_{n\to\infty}\sup_{t\in S_\omega}\|\bar u_{\eps_n,t}(\omega)\|_{L^1}
\leq\varliminf_{n\to\infty}\frac{1}{n}\sum_{k=1}^n
\sup_{t\in S_\omega}\|u_{\eps_k,t}(\omega)\|_{L^1},
\de
which implies  (\ref{Es77}) by (\ref{Es09}).
\end{proof}

Next, we study the weak continuity of generalized solutions.
We need the following technical lemma.
\bl\label{Le6}
Let $v\in C_w([0,T];L^2(\mR^d))$ (resp. $v\in C([0,T];L^2(\mR^d))$). If for some $R_N\to\infty$,
\be
\lim_{N\to\infty}\mathrm{ess.}\sup_{t\in[0,T]}\int_{|x|\geq R_N}|v_t|=0,\label{Cp5}
\ee
then $v\in C_w([0,T];L^1(\mR^d))$ (resp. $v\in C([0,T];L^1(\mR^d))$).
\el
\begin{proof}
We first prove that
\be
\mathrm{ess.}\sup_{t\in[0,T]}\int_{|x|\geq R_N}|v_t|=
\sup_{t\in[0,T]}\int_{|x|\geq R_N}|v_t|.\label{Cp6}
\ee
Let $\cS\subset[0,T]$ with full measure such that
$$
\sup_{t\in\cS}\int_{|x|\geq R_N}|v_t|=\mathrm{ess.}\sup_{t\in[0,T]}\int_{|x|\geq R_N}|v_t|.
$$
For $t\notin\cS$, let $\{t_k,k\in\mN\}\subset\cS$ converge to $t$. Since
$v_{t_k}$ weakly converges to $v_t$ in $L^2(\mR^d)$, by Banach-Saks theorem (cf. \cite{Du-Sc}),
there exists a subsequence still denoted by $t_k$ such that its Ces\`aro mean
$\bar v_{t_n}:=\frac{1}{n}\sum_{i=1}^nv_{t_k}$ strongly converges to $v_t$ in $L^2(\mR^d)$.
Thus, by Fatou's lemma,
$$
\int_{|x|\geq R_N}|v_t|\leq\lim_{n\to\infty}\int_{|x|\geq R_N}|\bar v_{t_n}|
\leq\sup_{t\in\cS}\int_{|x|\geq R_N}|v_t|,
$$
which then leads to (\ref{Cp6}).

Let $v\in C_w([0,T];L^2(\mR^d))$ and $\phi\in L^\infty(\mR^d)$. For $t_n\to t$,
we write
$$
\int (v_{t_n}-v_t)\phi=\int_{|x|\leq R_N} (v_{t_n}-v_t)\phi+\int_{|x|>R_N} (v_{t_n}-v_t)\phi.
$$
By (\ref{Cp5}) and (\ref{Cp6}), the second term can be arbitrarily small uniformly in $n$
for $N$ large enough. For fixed $N$, the first term goes to zero as $n\to\infty$.
The desired continuity then follows. If $v\in C([0,T];L^2(\mR^d))$ and $t_n\to t$, then
$$
\int |v_{t_n}-v_t|=\int_{|x|\leq R_N} |v_{t_n}-v_t|+\int_{|x|>R_N} |v_{t_n}-v_t|.
$$
As above, we have $v\in C([0,T];L^1(\mR^d))$.
\end{proof}
Using this lemma, we can prove the following result about the weak continuity of
generalized solutions. Our proof is adapted from \cite[p.206, Theorem 3]{Ro}.
\bp\label{Pro7}
Let $u\in L^2(\Omega; L^\infty(0,T;L^2(\mR^d)))$ be a generalized solution of SPDE (\ref{SPDE}).
Then there exists a version $\tilde u\in L^2(\Omega; C_w([0,T];L^2(\mR^d)))$ so that
$u_t(\omega,x)=\tilde u_t(\omega,x)$ $(\dif t\times P\times\dif x)$-a.s..
Moreover, if $u$ also satisfies
\be
\lim_{R\to\infty}\mE\left(\mathrm{ess.}\sup_{t\in[0,T]}\int_{|x|\geq R}|u_t|\right)=0,\label{Cp3}
\ee
then $\tilde u$ also belongs to
$L^1(\Omega; C_w([0,T];L^1(\mR^d)))$.
\ep
\begin{proof}
Let $\sD=\{\phi_1,\cdots,\phi_n,\cdots\}
\subset C^\infty_0(\mR^d)$ be a countable and dense subset of $L^2(\mR^d)$.
For each $\phi\in\sD$, we write the right hand side of (\ref{Gen}) as $\Phi_t(\phi)$.
Then $t\mapsto \Phi_t(\phi)$ is a continuous process and for
$(\dif t\times P)-$almost all $(t,\omega)\in[0,T]\times\Omega$
$$
\Phi_t(\phi)(\omega)=\int u_t(\omega)\phi.
$$
Let $\{r_1,r_2,\cdots,r_n\}$ be $n$ rational numbers. Then
\be
|r_i\Phi_t(\phi_i)(\omega)|\leq\|u_t(\omega)\|_{L^2}\|r_i\phi_i\|_{L^2}
\leq\mathrm{ess}\sup_{t\in[0,T]}\|u_t(\omega)\|_{L^2}\|r_i\phi_i\|_{L^2}.\label{Cp1}
\ee
Let $\sR$ be the collection of all finite many rational numbers $Q=\{r_1,r_2,\cdots,r_n\}$.
By the countability of $\sD$ and $\sR$ as well as the continuity of the left hand side,
there is a common null set $N$ such that for all $\omega\notin N$ and
all $t\in[0,T]$, $\phi\in\sD$, $Q\in\sR$, inequality (\ref{Cp1}) holds true.

Below, we fix such an $\omega\notin N$.
Let $\cL(\sD)$ be the linear space spanned by $\sD$. By the continuous dependence of
both sides of (\ref{Cp1}) in $Q\in\sR$, one can define a linear functional
$\hat\Phi_t$ on $\cL(\sD)$ such that
$$
\hat\Phi_t(\phi)(\omega)=\Phi_t(\phi)(\omega),\ \ \forall \phi\in\sD
$$
and
$$
|\hat\Phi_t(\phi)(\omega)|\leq\mathrm{ess}\sup_{t\in[0,T]}\|u_t(\omega)\|_{L^2}\|\phi\|_{L^2},
\ \ \forall \phi\in\cL(\sD).
$$
By Hahn-Banach theorem (cf. \cite{Du-Sc}), there exists a linear functional $\tilde\Phi_t$ such that
$$
\tilde\Phi_t(\phi)(\omega)=\hat\Phi_t(\phi)(\omega),\ \ \forall \phi\in\cL(\sD)
$$
and
\be
|\tilde\Phi_t(\phi)(\omega)|\leq\mathrm{ess}\sup_{t\in[0,T]}\|u_t(\omega)\|_{L^2}\|\phi\|_{L^2},
\ \ \forall \phi\in L^2(\mR^d).\label{Ep5}
\ee
By Riesz theorem, there exists a unique $\tilde u_t\in L^2(\mR^d)$ such that
$$
\tilde\Phi_t(\phi)(\omega)=\int\tilde u_t(\omega)\phi\ \mbox{ and }\ \|\tilde u_t\|_{L^2}
\leq \mathrm{ess}\sup_{t\in[0,T]}\|u_t(\omega)\|_{L^2}.
$$
Since for any $\phi\in\sD$, $t\mapsto\int\tilde u_t(\omega)\phi=
\tilde\Phi_t(\phi)(\omega)=\Phi_t(\phi)(\omega)$
is continuous, by (\ref{Ep5}), we also have for any $\phi\in L^2(\mR^d)$,
$$
t\mapsto\int\tilde u_t(\omega)\phi\mbox{ is continuous. }
$$
The first conclusion is then proven. The second conclusion follows from Lemma \ref{Le6}.
\end{proof}

Below, we give sufficient conditions for (\ref{Cp3}).
\bp\label{Pro8}
In anyone situation of Propositions \ref{Pro4}, \ref{Pro6}
and \ref{Pro5}, we also assume that
$u_0\in L^1(\Omega,\cF_0; L^1(\mR^d))$, $f\in L^1([0,T]\times\Omega\times\mR^d)$
and one of (\ref{C3}) and (\ref{C33}) hold. If
$$
u\in L^2(\Omega; L^\infty(0,T;L^2(\mR^d)))
\cap L^1(\Omega; L^\infty(0,T;L^1(\mR^d)))
$$
is a generalized solution of SPDE (\ref{SPDE}),
then (\ref{Cp3}) holds.
\ep
\begin{proof}We only consider the situation of Proposition \ref{Pro4}.
Let $\lambda_R(x)=\lambda(x/R)$, where $\lambda$ is a non-negative smooth function on $\mR^d$
with $\lambda(x)=1$ for $|x|\geq 2$ and $\lambda(x)=0$ for $|x|\leq 1$. Let $\chi_n(x)=\chi(x/n)$
be a cutoff function. Following the proof of (II) in Proposition \ref{Pro4},
we choose $\phi=\phi^R_n=\lambda_R\cdot\chi_n$ in (\ref{Cp2}). Then
\ce
\mE\left(\mathrm{ess.}\sup_{s\in[0,t]}\int|u_s|\phi^R_n\right)&\leq&
C\mE\int|u_0|\phi^R_n+C\mE\int^T_0\!\!\!\int |f|\phi^R_n\\
&&+C\mE\int^T_0\!\!\!\int|u|(|\p_j(\p_i\phi^R_n a^{ij})|+|\p_i\phi^R_n b^i|)\no\\
&&+C\mE\left(\int^T_0 \left|\int|u|~|\sigma^{il}\p_i\phi^R_n\right|^2 \right)^{1/2}.
\de
Notice that
\ce
|\p_i\phi_n^R|\leq\frac{C\chi_n 1_{|x|\geq R}}{1+|x|}
+\frac{C\lambda_R 1_{|x|\geq n}}{1+|x|}
\de
and
\ce
|\p_i\p_j\phi_n^R|\leq\frac{C\chi_n 1_{|x|\geq R}}{1+|x|^2}
+\frac{C\lambda_R 1_{|x|\geq n}}{1+|x|^2}+
\frac{C 1_{|x|\geq n} 1_{|x|\geq R}}{1+|x|^2}.
\de
Firstly letting $n\to\infty$ and then $R\to\infty$, as in the proof of Proposition \ref{Pro4}, we get
$$
\lim_{R\to\infty}\mE\left(\mathrm{ess.}\sup_{s\in[0,t]}\int_{|x|\geq 2R}|u_s|\right)
\leq\lim_{R\to\infty}\lim_{n\to\infty}
\mE\left(\mathrm{ess.}\sup_{s\in[0,t]}\int|u_s|\phi^R_n\right)=0.
$$
The proof is complete.
\end{proof}

\br
Using Lemma \ref{Le6} and Propositions \ref{Pro1}, \ref{Pro7} and \ref{Pro8}, we can improve
\cite[p.204, Corollary 1]{Ro} so that $u\in L^1(\Omega; C([0,T];L^1(\mR^d)))$ since all
the coefficients therein are bounded and $u\in L^2(\Omega; C([0,T];L^2(\mR^d)))$.
\er

\section{Application to Nonlinear Filtering}

Let $(\hat B_t)_{t\in[0,T]}$ and $(\tilde B_t)_{t\in[0,T]}$ be two independent
$d$ and $d_1$ dimensional standard Brownian motions on a standard filtered probability
space $(\Omega,\cF,P;(\cF_t)_{t\in[0,T]})$. Let $x_t$ denote the $d$-dimensional
unobservable signal and $y_t$ the $d_1$-dimensional observable signal. We assume that
$z_t=(x_t,y_t)$ obeys the following It\^o SDE:
$$
\dif\left(
\begin{array}{c}
x_t\\
y_t
\end{array}
\right)
=\left(
\begin{array}{c}
\hat b_t(z_t)\\
\tilde b_t(z_t)
\end{array}
\right)\dif t+
\left(
\begin{array}{cc}
\hat\sigma_t(z_t)& 0\\
0&\tilde\sigma_t(y_t)
\end{array}
\right)
\dif\left(
\begin{array}{c}
\hat B_t\\
\tilde B_t
\end{array}
\right),
$$
where $z_0=(x_0,y_0)$ is an $\cF_0$-measurable random variable and
the coefficients satisfy the following conditions:
\begin{enumerate}[{\bf (H1)}]
\item The regular conditional distribution of $x_0$ with respect to the
$\sigma$-algebra generated by $y_0$ is absolutely continuous with respect
to the Lebesgue measure on $\mR^d$ and the desity $\pi_0\in L^2(\Omega; L^2(\mR^d))$.
\item The functions $\hat b,\hat\sigma, \tilde b,\tilde\sigma$ satisfy the Lipschitz
conditions with respect to $z$ with constant $K$. Moreover, $\hat\sigma_t(x,y)$ is continuously
differentiable with respect to $x$ (not $z$) and its first derivatives with respect to $x^i$
satisfy the Lipschitz condition with respect to $x$ (not $z$) with constant $K$ independent of $y$.
\item $\tilde \sigma$ is non-singular and $\hat\sigma, \hat b(0,\cdot), \tilde\sigma,\tilde\sigma^{-1}, \tilde b$
are bounded by $K$.
\end{enumerate}
These assumptions will be forced throughout this section.

Let $\cF^y_t$ be the $P$-complete $\sigma$-algebra generated by $\{y_s,s\leq t\}$, which represents
the observation information. We want to get the conditional distribution of $x_t$
under $\cF^y_t$, i.e, to calculate
$$
\Pi_t(\omega,\Gamma):=P(x_t\in\Gamma|\cF^y_t),
$$
which is called the problem of filtering.

Let
\ce
&a^{ij}_t(\omega,x):=\hat\sigma^{ik}_t(x,y_t(\omega))\hat\sigma^{jk}_t(x,y_t(\omega))/2,\\
&h_t(\omega,x):=\tilde\sigma_t^{-1}(y_t(\omega))\tilde b_t(x,y_t(\omega)).
\de
We introduce the differential operators $\sL_t(\omega,x)$ and $\sM_t(\omega,x)$ by
\ce
\sL_t(\omega,x) u&:=&\p_i\p_j(a^{ij}_t(\omega,x)u)
-\p_i(\hat b^i_t(x,y_t(\omega))u)\\
&=&\p_i(a^{ij}_t(\omega,x)\p_ju)-\p_i(b^i_t(\omega,x)u),
\de
where $b^i_t(\omega,x):=\hat b^i_t(x, y_t(\omega))-\p_ja^{ij}_t(\omega,x)$, and
$$
\sM_t(\omega,x)u:=h_t(\omega,x) u.
$$

Define
$$
\rho_t:=\exp\left\{\int^t_0 h^k_s(x_s)\dif\tilde B^k_s+\frac{1}{2}\int^t_0|h^k_s(x_s)|^2\dif s\right\}
$$
and
$$
\bar B^k_t:=\tilde B^k_t+\int^t_0h^k_s(x_s)\dif s.
$$
By Girsanov's theorem, under the new probability measure
$$
\bar P(\dif\omega):=\rho^{-1}_t(\omega) P(\dif\omega),
$$
$\bar B_t$ is still a $d_1$-dimensional standard Brownian motion and
independent of $\hat B_t$. Moreover,
$$
\dif\left(
\begin{array}{c}
x_t\\
y_t
\end{array}
\right)
=\left(
\begin{array}{c}
\hat b_t(z_t)\\
0
\end{array}
\right)\dif t+
\left(
\begin{array}{cc}
\hat\sigma_t(z_t)& 0\\
0&\tilde\sigma_t(y_t)
\end{array}
\right)
\dif\left(
\begin{array}{c}
\hat B_t\\
\bar B_t
\end{array}
\right).
$$
The following lemma is taken from \cite[p. 228, Lemma 1.4]{Ro}.
\bl
Let $\bar\cF_t$ be the $\sigma$-algebra generated by $\{\bar B_s,s\leq t\}$. Then
$$
\cF^y_t=\bar\cF_t\vee\cF^y_0.
$$
\el

From this lemma, we know that $\bar B_t$ is a $d_1$-dimensional standard Brownian motion on
filtered probability space $(\Omega,\cF,\bar P;(\cF^y_t)_{t\in[0,T]})$. Moreover, it is clear that
the coefficients in $\sL$ and $\sM$ are measurable and $\cF^y_t$-adapted.
Consider the following SPDE:
$$
\dif u_t=\sL_t u_t\dif t+\sM^k_t u_t\dif\bar B_t^k,\ \ u_0=\pi_0.
$$
Under {\bf (H1)-(H3)}, by Theorem \ref{Main2} , there exists a unique non-negative generalized solution
in the class that
$$
u\in L^2(\Omega; C_w([0,T];L^2(\mR^d)))
\cap L^1(\Omega; C_w([0,T];L^1(\mR^d)))
$$
and
$$
\mE\left(\int^T_0\!\!\!\int|\hat\sigma^{il}\p_i u|^2\right)<+\infty.
$$

 We now give a representation for $u_t$.
\bp
For any $\phi\in L^\infty(\mR^d)$ and $t\in[0,T]$, we have
\be
\int u_t\phi=\mE^{\bar P}(\phi(x_t)\rho_t|\cF^y_t),\ \ P-a.s..\label{Rep}
\ee
\ep
\begin{proof}
By suitable approximation, we only need to prove (\ref{Rep}) for $\phi\in C^\infty_0(\mR^d)$.
We sketch the proof. As in Section 3, we define
$$
\hat\sigma^{ik}_{t,\eps}(\omega,\cdot):=(\sigma^{ik}_t(\cdot,y_t(\omega))*\rho_\eps)\chi_\eps,\
h^l_\eps:=(h^l*\rho_\eps)\chi_\eps
$$
and
$$
\hat b^i_{t,\eps}(\omega,\cdot):=[(\hat b^i_t(\cdot,y_t(\omega))
\wedge(1/\eps))\vee(-1/\eps)]*\rho_\eps,
$$
and consider the corresponding approximation equation:
\be
\dif u_{\eps,t}=\sL_{t,\eps} u_{\eps,t}\dif t+\sM^k_{t,\eps} u_{\eps,t}\dif\bar B_t^k,
\ \ u_{\eps,0}=\pi_0.\label{P4}
\ee
By \cite[p.203, Theorem 1]{Ro} (see also \cite{Ku}),
the unique solution of equation (\ref{P4}) can be represented by
\be
\int u_{\eps,t}\phi=\mE^{\bar P}(\phi(x_{\eps,t})\rho_{\eps,t}|\cF^y_t),\label{Rep1}
\ee
where $x_{\eps,t}$ solves the following SDE:
\ce
x_{\eps,t}=x_0+\int^t_0\hat b_{s,\eps}(x_{\eps,s})\dif s
+\int^t_0\hat\sigma_{s,\eps}(x_{\eps,s})\dif\hat B_s
\de
and
$$
\rho_{\eps,t}=1+\int^t_0 \rho_{\eps,s} h^k_{s,\eps}(x_{\eps,s})\dif \bar B^k_s,
$$
i.e.,
$$
\rho_{\eps,t}=\exp\left\{\int^t_0 h^k_{s,\eps}(x_{\eps,s})\dif \bar B^k_s
-\frac{1}{2}\int^t_0|h^k_{s,\eps}(x_{\eps,s})|^2\dif s\right\}.
$$
It is now standard to prove that
$$
\sup_{\eps\in(0,1)}\mE\left(\sup_{s\in[0,T]}|x_{\eps,s}|^2\right)<+\infty.
$$
Using this estimate, we can prove that for any $\delta>0$
$$
\lim_{\eps\to0}P\left(\sup_{t\in[0,T]}|x_{\eps,t}-x_t|\geq\delta\right)=0
$$
and
$$
\lim_{\eps\to0}P\left(\sup_{t\in[0,T]}|\rho_{\eps,t}-\rho_t|\geq\delta\right)=0,
$$
where we have used that
$$
\rho_t=\exp\left\{\int^t_0 h_s(x_s)\dif\bar B_s-\frac{1}{2}\int^t_0|h^l_s(x_s)|^2\dif s\right\}.
$$
On the other hand, as in the proof of Proposition \ref{Th1}, one knows that
$$
u_{\eps,t}\stackrel{\eps\to 0}\longrightarrow u_t
\mbox{ weakly in $L^2([0,T]\times\Omega\times\mR^d)$}.
$$
Now taking weak limits for (\ref{Rep1}), we obtain
$$
\int u_t\phi=\mE^{\bar P}(\phi(x_t)\rho_t|\cF^y_t),\ \ (\dif t\times P)-a.s.
$$
Since the left hand side is continuous and the right hand side also admits a
continuous version (cf. \cite[p.206, Theorem 3]{Ro}), representation (\ref{Rep}) now follows.
\end{proof}
Our main result in this section is:
\bt
Under {\bf (H1)}-{\bf (H3)}, the conditional distribution $\Pi_t(\omega,\Gamma)$ has a density
$\pi_t(\omega,\cdot)\in C_w([0,T];L^1(\mR^d))$ with respect to the Lebesgue measure almost surely. It is given by
\be
\pi_t(\omega,x)=\frac{u_t(\omega,x)}{\int u_t(\omega,x)\dif x}.\label{For3}
\ee
Moreover, for any $\phi\in C^\infty_0(\mR^d)$, $\pi_t(\phi)=\int \phi\pi_t$
satisfies the following non-linear SPDE:
\be
\pi_t(\phi)=\pi_0(\phi)+\int^t_0\pi_s(\sL^*_s\phi)\dif s+
\int^t_0\big[\pi_s(\sM^{k*}_s\phi)-\pi_s(h^k_s)\pi_s(\phi)\big]\dif\check B^k_s,\label{Eq4}
\ee
where $\dif\check B^k_t=\dif\bar B^k_t-\pi_t(h^k_t)\dif t$ and $\sL^*_s$ and
$\sM^{k*}_s$ are their respective adjoint operators.
\et
\begin{proof}
By (\ref{Rep}) and Bayes' formula about the conditional expectations
(cf. \cite[p.224, Theorem 1]{Ro}), we  have
$$
\mE^P(\phi(x_t)|\cF^y_t)=\frac{\mE^{\bar P}(\phi(x_t)\rho_t|\cF^y_t)}
{\mE^{\bar P}(\rho_t|\cF^y_t)}=\left(\int u_t\right)^{-1}\int u_t\phi.
$$
Formula (\ref{For3}) follows.

Observe that
$$
\rho_t=1+\int^t_0 h^k_s(x_{s})\rho_s \dif \bar B^k_s.
$$
Taking the conditional expectation with respect to $\cF^y_t$, we get
\ce
\mE^{\bar P}(\rho_t|\cF^y_t)&=&1+\mE^{\bar P}\left(\int^t_0 h^k_s(x_{s})
\rho_s \dif \bar B^k_s\Big|\cF^y_t\right)\\
&=&1+\int^t_0 \mE^{\bar P}(h^k_s(x_{s})
\rho_s|\cF^y_s) \dif \bar B^k_s\\
&=&1+\int^t_0 \mE^{P}(h^k_s(x_{s})|\cF^y_s)\mE^{\bar P}(\rho_s|\cF^y_s) \dif \bar B^k_s,
\de
where the second equality is due to the property of stochastic integrals and the
third equality is due to the Bayes' formula.

In view of $h^k_s(x_{s})=\tilde\sigma_s^{-1}(y_s)\tilde b_s(x_s,y_s)$, by certain approximation,
we have
$$
\mE^{P}(h^k_s(x_{s})|\cF^y_s)=\pi_s(h^k_s).
$$
Hence,
$$
\mE^{\bar P}(\rho_t|\cF^y_t)=\exp\left\{\int^t_0\pi_s(h^k_s)\dif \bar B^k_s
-\frac{1}{2}\int^t_0|\pi_s(h^k_s)|^2\dif s\right\}.
$$
Since $\int u_t=\mE^{\bar P}(\rho_t|\cF^y_t)$, equation (\ref{Eq4}) now follows
by It\^o's formula.
\end{proof}

\br
We can also consider the filtering problem in the cases of Theorems \ref{Main1}, \ref{Main3} and
\ref{Main4}. In particular, in the case of Theorem \ref{Main3}, we can even allow some
singularity of $\hat b$ in $x$.
\er

\section{A Degenerate nonlinear SPDE}

Let $a^{ij}$ be given by
$$
a^{ij}_t(\omega,x)=\hat\sigma^{il}_t(\omega,x)\hat\sigma^{jl}_t(\omega,x)
$$
and
$$
\sL_t(\omega)u:=\p_i(a^{ij}_t(\omega,x)\p_ju)
+\p_i(b^i_t(\omega,x)u),\ \
\sM^l_t(\omega)u:=\sigma^{il}_t(\omega,x)\p_iu.
$$
In this section, we consider the following SPDE
\be
\dif u_t=(\sL_tu_t+f_t(u_t))\dif t+(\sM^l_tu_t+g^l_t(u_t))\dif B^l_t, \
u_0(\omega,x)=\varphi(\omega,x),\label{NonLinear}
\ee
where
$$
f: [0,T]\times\Omega\times\mR^d\times\mR\to\mR,\
g: [0,T]\times\Omega\times\mR^d\times\mR\to l^2.
$$
are $\cM\times\cB(\mR^d\times\mR)$-measurable functions.

Our main result in this section is:
\bt\label{Th2}
Assume that for some $\a>1/2$
$$
|\hat\sigma^{il}_t\xi_i|^2\geq\a|\sigma^{il}_t\xi_i|^2
$$
and the following conditions hold: for some $q>1$
$$
\left\{
\begin{aligned}
&\frac{\hat\sigma^{i\cdot}}{1+|x|}, \p_i\hat\sigma^{i\cdot}
\in L^2(0,T; L^\infty(\Omega\times\mR^d; l^2)),\\
&\frac{|b^i|}{1+|x|}, \div b\in L^1(0,T; L^\infty(\Omega\times\mR^d)),\\
&\p_kb^i\in L^1(0,T; L^2(\Omega;L^{2}_{loc}(\mR^d))),\\
&\sigma^{i\cdot}, \p_i\sigma^{i\cdot}
\in L^{2q}(0,T; L^\infty(\Omega\times\mR^d; l^2)),
\end{aligned}
\right.
$$
and for some $K>0$ and $\gamma\in L^2([0,T]\times\Omega\times\mR^d)$
\be
&|f_t(\omega,x,z)-f_t(\omega,x,z')|+\|g_t(\omega,x,z)-g_t(\omega,x,z')\|_{l^2}
\leq K |z-z'|,\label{Con3}\\
&|f_t(\omega,x,z)|+\|g_t(\omega,x,z)\|_{l^2}
\leq K|z|+\gamma_t(\omega,x).\label{Con4}
\ee
Then for any $u_0\in L^2(\Omega,\cF_0; L^2(\mR^d))$, there exists a unique generalized solution with
\be
u\in L^2(\Omega; C_w([0,T]; L^2(\mR^d)))\label{Con1}
\ee
and
\be
\hat\sigma^{i\cdot}\p_i u, \sigma^{i\cdot}\p_i u
\in L^2([0,T]\times\Omega; L^2(\mR^d;l^2)).\label{Con2}
\ee
\et
\begin{proof}
(Uniqueness): The uniqueness is a conclusion of the maximum principle. In fact, let
$u$ and $\tilde u$ be two generalized solutions of nonlinear SPDE (\ref{NonLinear})
with the same initial values and satisfy (\ref{Con1}) and (\ref{Con2}).
It is easy to see that
$$
v:=u-\tilde u
$$
satisfies the following linear equation:
$$
\dif v_t=(\sL_tv_t+c_tv_t)\dif t+(\sM^l_tv_t+h^l_tv_t)\dif B^l_t, \ \
v_0(\omega,x)=0 ,
$$
where
$$
c_t(\omega,x)=\int^1_0(\p_z f_t)(\omega,x, \theta(u-\tilde u)+\tilde u)\dif \theta
$$
and
$$
h^l_t(\omega,x)=\int^1_0(\p_z g^l_t)
(\omega,x, \theta(u-\tilde u)+\tilde u)\dif \theta.
$$
By assumption (\ref{Con3}), we know that
$$
c\in L^\infty([0,T]\times\Omega\times\mR^d)),\ \
h\in L^\infty([0,T]\times\Omega\times\mR^d; l^2)).
$$
Hence, by Proposition \ref{Pro5}, we have $v\equiv 0$. The uniqueness then follows.

(Existence): We use the Picard iteration method and a priori estimate (\ref{Ep1}).
Let $u^0_t(\omega,x)=\varphi(\omega,x)$. Consider the following approximation equation:
\be
\dif u^{n}_t=(\sL_tu^n_t+f_t(u^{n-1}_t))\dif t+(\sM^l_tu^n_t+g^l_t(u^{n-1}_t))\dif B^l_t, \ \
u^n_0(\omega,x)=\varphi(\omega,x).\label{NonLinearApp}
\ee
By (\ref{PL1}) and (\ref{Con4}) we have
\ce
&&\mE\left(\sup_{s\in[0,t]}\int|u^n_s|^2\right)
+\mE\left(\int^t_0\!\!\!\int(|\hat\sigma^{il}\p_i u^n|^2+|\sigma^{il}\p_i u^n|^2)\right)\\
&&\qquad\leq C\mE\int|u_0|^2+C\mE\int^T_0\!\!\!\int|f_s(u^{n-1}_s)|^2
+C\mE\int^T_0\!\!\!\int|g^l_s(u^{n-1}_s)|^2\\
&&\qquad\leq C\mE\int|u_0|^2+C\int^T_0\ell_s\dif s+C\mE\int^t_0\!\!\!\int|u^{n-1}_s|^2,
\de
By Gronwall's inequality, we get the following uniform estimates:
\be
\mE\left(\sup_{s\in[0,T]}\int|u^n_s|^2\right)
+\mE\left(\int^T_0\!\!\!\int(|\hat\sigma^{il}\p_i u^n|^2+|\sigma^{il}\p_i u^n|^2)\right)
\leq C,\label{Ep2}
\ee
where $C$ is independent of $n$.

Set now
$$
v^{n,m}:=u^n-u^m.
$$
Then, by (\ref{PL1}) again, we have
\ce
\mE\left(\mathrm{ess.}\sup_{s\in[0,t]}\int |v^{n,m}_s|^2\right)&\leq&
C\mE\int^t_0\!\!\!\int|f_s(u^{n-1}_s)-f_s(u^{m-1}_s)|^2\\
&&+C\mE\int^t_0\!\!\!\int|g^l_s(u^{n-1}_s)-g^l_s(u^{m-1}_s)|^2.\\
&\leq& C\int^t_0 \mE\int|v^{n-1,m-1}_s|^2.
\de
Set
$$
\Phi_t:=\varlimsup_{n,m\to\infty}\mE\left(\mathrm{ess.}\sup_{s\in[0,t]}\int |v^{n,m}_s|^2\right).
$$
Then by (\ref{Ep2}) and Fatou's lemma, we have
$$
\Phi_t\leq C\int^t_0 \Phi_s\dif s,
$$
which implies that
$$
\varlimsup_{n,m\to\infty}\mE\left(\mathrm{ess.}\sup_{s\in[0,T]}\int |v^{n,m}_s|^2\right)
=\Phi_T=0.
$$
So, there is a $u\in L^2(\Omega; L^\infty(0,T; L^2(\mR^d)))$ such that
$$
\lim_{n\to\infty}\mE\left(\mathrm{ess.}\sup_{s\in[0,T]}\int |u^{n}_s-u_s|^2\right)=0.
$$
By passing to the limits for (\ref{NonLinearApp}), we obtain that $u$ is a generalized
solution. (\ref{Con1}) is due to Proposition \ref{Pro7}.
Estimate (\ref{Con2}) follows from (\ref{Ep2}).
\end{proof}

\br
If $g_t(\omega,x,0)=0$, $f_\cdot(\cdot,\cdot,0)\in L^1([0,T]\times\Omega\times\mR^d)$
and $u_0\in L^1(\Omega\times\mR^d)$, then the unique solution in Theorem \ref{Th2}
also belongs to $L^1(\Omega; C_w([0,T]; L^1(\mR^d)))$.
\er

\section{Appendix: Proof of Lemma \ref{Le2}}

We only prove (\ref{Lm1}). For $R>0$, let $B_R:=\{x\in\mR^d: |x|\leq R\}$. Below, we simply write
$$
\|f\|_{L^{r_1}(0,T; L^{r_2}(\Omega, L^{r_3}(B_R)))}=:\|f\|_{r_1,r_2,r_3;R}
$$
and
\ce
\mS_1&:=&L^{q_1}(0,T; L^{q_2}(\Omega, W^{1,q_3}_{loc}(\mR^d))),\\
\mS_2&:=&L^{p_1}(0,T; L^{p_2}(\Omega, L^{p_3}_{loc}(\mR^d))),\\
\mS_3&:=&L^{r_1}(0,T; L^{r_2}(\Omega, L^{r_3}_{loc}(\mR^d))).
\de

Notice that
$$
[\rho_\eps,b^i\p_i](u)=\int (b^i(y)-b^i(x))u(y)\p_i\rho_\eps(x-y)\dif y
-\int \div b(y) u(y)\rho_\eps(x-y)\dif y.
$$
If $b$ and $u$ are smooth in $x$, then it is easy to see that for every $x\in\mR^d$,
$$
\int (b^i(y)-b^i(x))u(y)\p_i\rho_\eps(x-y)\dif y
\stackrel{\eps\to0}{\longrightarrow} \div b(x) u(x),
$$
which implies by the dominated convergence theorem,
$$
\lim_{\eps\to 0}\|[\rho_\eps,b^i\p_i](u)\|_{r_1,r_2,r_3;R}=0.
$$

{\bf (Case: $p_1,q_1,p_2,q_2,p_3,q_3<+\infty$)}. It is enough to show that
\be
\mS_1\times\mS_2\ni (b,u)\mapsto [\rho_\eps,b^i\p_i](u)\in \mS_3\label{Ep7}
\ee
is uniformly continuous with respect to $\eps$.
First of all, by H\"older's inequality, we have for any $R>0$,
\be
\left\|\int \div b(y) u(y)\rho_\eps(x-y)\dif y\right\|_{r_1,r_2,r_3;R}
\leq \|\div b\|_{q_1,q_2,q_3;R+1}\|u\|_{p_1,p_2,p_3;R+1}.\label{Ep6}
\ee
Observing that
$$
|b^i(y)-b^i(x)|\leq |y-x|\int^1_0|\nabla b^i|(y+\theta(x-y))\dif\theta
$$
and
$$
\eps|\nabla\rho_\eps|(x)\leq C\rho_{2\eps}(x),
$$
where $C$ is independent of $\eps$, we have
\ce
\gamma_\eps(x)&:=&\left|\int (b^i(y)-b^i(x))u(y)\p_i\rho_\eps(x-y)\dif y\right|\\
&\leq&C\int|u(y)|\int^1_0|\nabla b|(y+\theta(x-y))\dif\theta\rho_{2\eps}(x-y)\dif y\\
&=&C\int|u(x-y)|\int^1_0|\nabla b|(x-y+\theta y)\dif\theta\rho_{2\eps}(y)\dif y.
\de
Hence, by H\"older's inequality again,
$$
\|\gamma_\eps\|_{r_1,r_2,r_3;R}\leq C\|u\|_{p_1,p_2,p_3;R+1}\|\nabla b\|_{q_1,q_2,q_3;R+1},
$$
which together with (\ref{Ep6}) yields (\ref{Ep7}).

{\bf (Case: any of $p_1,q_1,p_2, q_2,p_3, q_3=+\infty$)}. Without loss of generality,
we assume $p_1=p_2=p_3=+\infty$ and $q_1, q_2, q_3<+\infty$.
In this case, let $u_\delta:=u*\rho_\delta$. It is enough to prove that
$$
\lim_{\delta\to 0}\sup_{\eps\in(0,1)}|[\rho_\eps,b^i\p_i](u-u_\delta)|=0.
$$
Since $\|u_\delta\|_{L^\infty(B_R)}\leq \|u\|_{L^\infty(B_{R+1})}$, and for almost all $x\in\mR^d$,
$$
u_\delta(x)\stackrel{\delta\to0}{\longrightarrow}u(x),
$$
by the dominated convergence theorem, one can see that
$$
\lim_{\delta\to 0}\sup_{\eps}\left\|\int \div b(y) (u(y)-u_\delta(y))
\rho_\eps(x-y)\dif y\right\|_{r_1,r_2,r_3;R}=0.
$$
Similarly,
\ce
&&\lim_{\delta\to 0}\sup_{\eps}\left\|\int (b^i(y)-b^i(x))(u(y)-u_\delta(y))
\p_i\rho_\eps(x-y)\dif y\right\|_{r_1,r_2,r_3;R}\\
&\leq&C\lim_{\delta\to 0}\sup_{\eps}\left\|\int|u-u_\delta|(x-y)
\int^1_0|\nabla b|(x-y+\theta y)\dif \theta
\rho_{2\eps}(y)\dif y\right\|_{r_1,r_2,r_3;R}=0.
\de
The proof of Lemma \ref{Le2} is thus complete.

\vspace{5mm}

{\bf Acknowledgements:}

The author would like to thank Professor Benjamin Goldys for
providing him an excellent environment to work in the University of New South Wales.
His work is supported by ARC Discovery grant DP0663153 of Australia and
NSF of China (No. 10871215).

\end{document}